\documentclass[preprint,3p]{article}

\usepackage[separate-uncertainty=true,multi-part-units=single,allow-number-unit-breaks]{siunitx}

\usepackage[margin=1in]{geometry}
\usepackage[sort&compress,square,numbers]{natbib}

\usepackage{upgreek}
\usepackage{textcomp}
\usepackage{stmaryrd}
\usepackage{verbatim} 
\usepackage{float, booktabs, bm, pgfplots, subcaption}
\usepackage{amsmath, amsfonts, amssymb, times}
\usepackage{graphicx}
\usepackage[vlined,linesnumbered,ruled]{algorithm2e}
\usepackage{todonotes}
\usepackage[capitalize]{cleveref}
\usepackage{multirow}

\usepgfplotslibrary{colorbrewer}
\pgfplotsset{
  log x ticks with fixed point/.style={
      xticklabel={
        \pgfkeys{/pgf/fpu=true}
        \pgfmathparse{exp(\tick)}
        \pgfkeys{/pgf/fpu=false}
      }
  },
  log y ticks with fixed point/.style={
      yticklabel={
        \pgfkeys{/pgf/fpu=true}
        \pgfmathparse{exp(\tick)}%
        \pgfmathprintnumber[fixed relative, precision=3]{\pgfmathresult}
        \pgfkeys{/pgf/fpu=false}
      }
  }
}

\pgfplotsset{width=8.0cm,
	     scaled y ticks=false,
             cycle list/Dark2,
             cycle multiindex* list={
	       mark list*\nextlist
	       linestyles\nextlist
	       Dark2\nextlist
	     },
	    }

\pgfplotsset{
  log x ticks with fixed point/.style={
      xticklabel={
        \pgfkeys{/pgf/fpu=true}
        \pgfmathparse{exp(\tick)}%
        \pgfmathprintnumber[fixed relative, precision=3]{\pgfmathresult}
        \pgfkeys{/pgf/fpu=false}
      }
  },
  log y ticks with fixed point/.style={
      yticklabel={
        \pgfkeys{/pgf/fpu=true}
        \pgfmathparse{exp(\tick)}%
        \pgfmathprintnumber[fixed relative, precision=3]{\pgfmathresult}
        \pgfkeys{/pgf/fpu=false}
      }
  }
}

\makeatletter
\renewcommand{\@algocf@capt@plain}{above}
\makeatother

\SetKwIF{If}{ElseIf}{Else}{if}{:}{else if}{else}{endif}
\SetKw{And}{\textbf{and}}

\newcommand{\T}{\mathrm{T}}

\newcommand{\p}{\mathrm{p}}

\newcommand{\comp}{\mathrm{c}}
\newcommand{\tens}{\mathrm{t}}

\newcommand{\tr}{\mathrm{tr}}
\newcommand{\etal}{\textit{et al}. }
\newcommand{\ie}{\textit{i}.\textit{e}. }
\newcommand{\eg}{\textit{e}.\textit{g}. }
\newcommand{\cf}{\textit{c}.\textit{f}. }

\newcommand{\norm}[1]{\left\lVert{#1}\right\rVert}

\newcommand{\mrm}{\mathrm}
\newcommand{\mbf}{\mathbf}
\newcommand{\bs}{\boldsymbol}

\newcommand{\fetwo}{FE$^\mathrm{2}$}
\newcommand{\strain}{\ensuremath{\bs{\varepsilon}}}
\newcommand{\stress}{\ensuremath{\bs{\sigma}}}
\newcommand{\disp}{\ensuremath{\mbf{u}^\omega}}

\newcommand{\data}{\ensuremath{\mathcal{D}}}
\newcommand{\weights}{\ensuremath{\mathbf{w}}}

\newcommand{\model}{\ensuremath{\mathcal{M}}}
\newcommand{\extractor}{\ensuremath{\mathcal{F}}}
\newcommand{\props}{\ensuremath{\boldsymbol{\theta}}}
\newcommand{\state}{\ensuremath{\boldsymbol{\alpha}}}
\newcommand{\features}{\ensuremath{\boldsymbol{\varphi}}}

\usepackage{authblk}

\title{Machine learning of evolving physics-based material models for multiscale solid mechanics}

\author[1]{I.B.C.M. Rocha}
\author[2]{P. Kerfriden}
\author[1]{F.P. van der Meer}
\affil[1]{Delft University of Technology, Faculty of Civil Engineering and Geosciences, P.O. Box 5048, 2600GA Delft, The Netherlands}
\affil[2]{Mines Paris, PSL University, Centre des mat\'{e}riaux, 63-65 Rue Henri-Auguste Desbrueres BP87, F-91003 \'{E}vry, France}
\date{}                     
\begin{document}

\maketitle

\begin{abstract}
In this work we present a hybrid physics-based and data-driven learning approach to construct surrogate models for concurrent multiscale simulations of complex material behavior. We start from robust but inflexible physics-based constitutive models and increase their expressivity by allowing a subset of their material parameters to change in time according to an evolution operator learned from data. This leads to a flexible hybrid model combining a data-driven encoder and a physics-based decoder. Apart from introducing physics-motivated bias to the resulting surrogate, the internal variables of the decoder act as a memory mechanism that allows path dependency to arise naturally. We demonstrate the capabilities of the approach by combining an FNN encoder with several plasticity decoders and training the model to reproduce the macroscopic behavior of fiber-reinforced composites. The hybrid models are able to provide reasonable predictions of unloading/reloading behavior while being trained exclusively on monotonic data. Furthermore, in contrast to traditional surrogates mapping strains to stresses, the specific architecture of the hybrid model allows for lossless dimensionality reduction and straightforward enforcement of frame invariance by using strain invariants as the feature space of the encoder.
\end{abstract}


\section{Introduction}
\label{SEintroduction}

Recent advances in materials science and manufacturing techniques are paving the way for the design of materials with highly-tailored microstructures, including metamaterials \cite{kumarInversedesignedSpinodoidMetamaterials2020,bessaBayesianMachineLearning2019}, novel composite material systems \cite{gantenbeinSpinPrintingLiquidCrystal2021,woigkFlaxbasedNaturalComposites2022}, printed cementitious materials \cite{xuUnderstandingDeformationFracture2022} and multifunctional living materials \cite{gantenbeinThreedimensionalPrintingMycelium2023a}. The common thread in these new developments is a shift from traditional design focused on tailoring structures to material constraints towards tailoring material microstructures to macroscopic constraints. This shift in turn requires the development of highly-detailed models of material behavior across spatial scales and a shift to virtual structural certification, as trial-and-error design becomes infeasible \cite{telgenTopologyOptimizationGraded2022,khajehtourianSoftAdaptiveMechanical2021,furtadoMethodologyGenerateDesign2021}.

Scale bridging has been traditionally performed through a bottom-up approach: physics-based constitutive models at smaller scales are calibrated using experiments and used to perform numerical simulations (using \eg the Finite Element (FE) method) on representative lower-scale domains from which higher-scale physics-based models can be calibrated \cite{vandermeerMicromechanicalValidationMesomodel2016,krauklisPredictionOrthotropicHygroscopic2019}. However, physics-based constitutive models come with \textit{a priori} assumptions that often fail to reproduce complex lower-scale behavior \cite{vandermeerMicromechanicalValidationMesomodel2016}. The alternative is to opt for an \fetwo\ (or Computational Homogenization) approach: lower-scale FE models are embedded at every Gauss point of a higher-scale model and material behavior is directly upscaled with no constitutive assumptions at the higher scale \cite{feyelMultiscaleFE2Elastoviscoplastic1999,kouznetsovaApproachMicromacroModeling2001,geersMultiscaleComputationalHomogenization2010}. Yet, the computational cost associated with repeatedly solving a large number of micromodels quickly becomes a bottleneck, in particular for many-query procedures such as design exploration and optimization that require several higher-scale simulations to be performed.

Since the bottleneck of \fetwo\ lies in computing lower-scale models, a popular approach to reduce computational effort is to substitute the original FE micromodels with either structure-preserving reduced-order models \cite{ferreiraAdaptivityClusteringbasedReducedorder2022,gouryAutomatisedSelectionLoad2016,ryckelynckPrioriHyperreductionMethod2005,ghavamianPODDEIMModel2017,rochaAdaptiveDomainbasedPOD2020,danielUncertaintyQuantificationIndustrial2022,scanffWeaklyinvasiveLATINPGDSolving2022} or purely data-driven surrogates \cite{ghaboussiKnowledgeBasedModeling1991,lefikArtificialNeuralNetworks2009,leComputationalHomogenizationNonlinear2015,bessaFrameworkDatadrivenAnalysis2017,rochaMicromechanicsbasedSurrogateModels2020,wangDeepLearningFramework2022} trained offline. More recently, Recurrent Neural Networks (RNN) have become the model of choice especially for strain path-dependent materials, with a large body of literature dedicated to their use and tuning to different applications \cite{ghavamianAcceleratingMultiscaleFinite2019,mozaffarDeepLearningPredicts2019,gorjiPotentialRecurrentNeural2020,abueiddaDeepLearningPlasticity2021,chenRecurrentNeuralNetworks2021,logarzoSmartConstitutiveLaws2021,borkowskiRecurrentNeuralNetworkbased2022}. RNNs can reproduce complex long-term time dependencies in material behavior by learning latent representations of the material state, making them fast and flexible surrogates. However, these learned representations are not a priori related to actual thermodynamic internal state variables and the model is therefore poorly interpretable (see \cite{koeppeExplainableArtificialIntelligence2021} for an interesting discussion on the subject). Furthermore, training for path dependency requires sampling from a potentially infinite-dimensional space of arbitrarily-long strain paths. This means training RNNs to reproduce complex material behavior often requires an inordinate amount of data (\textit{curse of dimensionality}) and their purely data-driven nature limits their ability to extrapolate away from paths seen during training.

In order to address these drawbacks, a growing number of recent works are shifting focus to models with a fusion of data-driven and physics-based components. Inspired by physics-informed neural networks (\cite{raissiPhysicsinformedNeuralNetworks2019}), the authors in \cite{masiThermodynamicsbasedArtificialNeural2021} opt for data-driven models with physics-inspired bias by enforcing thermodynamic principles in a weak sense through an augmented loss function. In a similar vein, the model in \cite{linkaConstitutiveArtificialNeural2021} learns hyperelasticity by linking together several carefully crafted neural nets to represent quantities with clear physical meaning, improving the interpretability of the resulting model. In \cite{vlassisSobolevTrainingThermodynamicinformed2021} the authors extend a similar hyperelastic surrogate with a network that learns plastic flow direction and the evolution of a yield surface parametrized by a level set function, resulting in a hyperelastic-plastic model with superior extrapolation capabilities. A common thread in the aforementioned approaches, however, is that their learning architectures are heavily dependent on the type of model being learned (\eg hyperelasticity, plasticity), making extensions to other models a convoluted task. In contrast, the authors in \cite{liuDeepMaterialNetwork2019,liuCellDivisionDeep2021} propose a surrogate for heterogeneous micromodels constructed by directly employing unmodified versions of the constitutive models used for the micro constituents and using a customized network architecture to infer a homogenization operator from data that combines their responses. Nevertheless, the method employs a highly-specialized iterative online prediction routine requiring extra implementation effort and with increased computational overhead when compared to that of traditional surrogates mapping strains to stresses. Finally, in \cite{wangCooperativeGameAutomated2019,flaschelUnsupervisedDiscoveryInterpretable2021,flaschelDiscoveringPlasticityModels2022} a dictionary of candidate physics-based models is assumed and the role of machine learning shifts instead to that of performing model selection and/or design of experiments.

\begin{figure}
\centering
\includegraphics[scale=1]{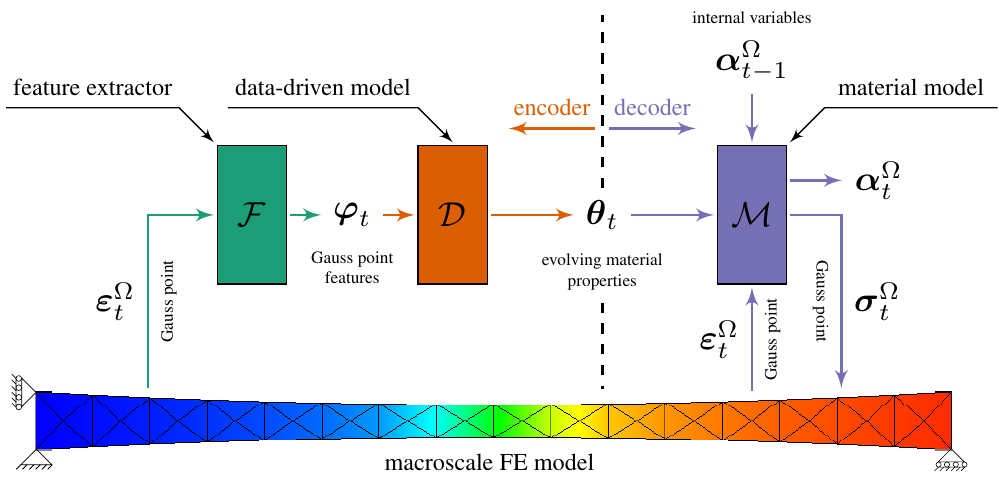}
\caption{The hybrid surrogate combining a data-driven encoder for material parameters and a physics-based material model decoder.}
\label{FIapproach}
\end{figure}

In this work we explore an alternative approach for constructing hybrid surrogate models for path-dependent multiscale simulations. We start from the premise that existing physics-based models --- \eg the ones used to describe microscale constituents --- are not flexible enough to reproduce macroscale behavior but nonetheless encapsulate crucial physical features such as frame invariance and loading/unloading conditions. It is our aim to avoid learning these features directly from data, as that would require either an excessively large dataset or a highly-specialized learning architecture. We therefore opt for keeping the constitutive model as intact as possible and instead increasing flexibility by allowing some (or all) of its material parameters to evolve in time. The resulting model can be seen in \cref{FIapproach}: a data-driven encoder that learns the evolution of a set of material properties is linked to a physics-based material model decoder that maps strains to stresses. In contrast to other strategies in literature, we keep the architecture as general as possible: a general feature extractor parses macroscopic strains into features for the encoder --- which can be as simple as the strains themselves or other derived quantities (\eg strain invariants) ---  and any type of constitutive model can in principle act as decoder (\eg hyperelasticity, plasticity, damage). By relegating stress computations to the decoder, we effectively introduce physics-based bias to the model.\footnote
{
In purely data-driven surrogates, we accept some bias in exchange for reduced variance --- \eg by employing regularization or adopting prior distributions for model parameters \cite{bishopPatternRecognitionMachine2006} --- in order to counter overfitting and improve generalization. But in that case the bias is merely a way to reduce complexity, with no physical interpretation and no \textit{a priori} impact on the extrapolation capabilities of the model.
}
Furthermore, by letting the material model handle the evolution of its own internal variables, the model benefits from a recurrent component with interpretable memory structure that allows path dependency to arise naturally. The strategy we explore here is related to the one we propose in \cite{maiaPhysicallyRecurrentNeural2022}, but in that work we let an encoder learn local strain distributions for several virtual material points with fixed properties. We see the two approaches as being complementary, and therefore with potential for being used in combination to form a flexible range of hybrid surrogates.


The remainder of the work is organized as follows. \cref{SEfe2} contains a primer on concurrent multiscale (\fetwo) modeling and discusses the difficulties of training purely data-driven surrogates. In \cref{SEmodel}, we particularize the model of \cref{FIapproach} to the case of a feedforward neural network encoder and discuss aspects related to offline training and online numerical stabilization. In \cref{SEexamples} we assess the performance of the hybrid model in reproducing the behavior of fiber-reinforced composites using different encoder features and decoder models. Finally, some concluding remarks and future research directions are discussed in \cref{SEconclusions}.

\section{Concurrent multiscale (\fetwo) modeling}
\label{SEfe2}

In this section we present a short discussion on \fetwo\ modeling. The goal is not to be comprehensive --- the interested reader is referred to \cite{kouznetsovaApproachMicromacroModeling2001,geersMultiscaleComputationalHomogenization2010} for detailed discussions on the subject --- but rather to expose the computational bottleneck associated with the method and pinpoint where surrogate models can be used to alleviate the issue. We then demonstrate how a Recurrent Neural Network (RNN) can be used as surrogate model and showcase some of the difficulties associated with their training and their extrapolation capabilities.

\subsection{Scale separation and coupling}
\label{SEfe2primer}

In \fetwo\ we assume the problem being solved can be split into a homogeneous macroscopic domain $\Omega$ and a heterogeneous microscopic domain $\omega\ll\Omega$ where small-scale geometric features are resolved. Here we opt for a first-order homogenization approach assuming the displacements on both scales can be related by: 
\begin{equation}
\disp = \strain^\Omega\mbf{x}^\omega + \widetilde{\mbf{u}}
\label{EQfe2disp}
\end{equation}
\noindent where microscopic displacements $\disp$ are split into a linear contribution proportional to the macroscopic strains $\strain^\Omega$ and a fluctuation term $\widetilde{\mbf{u}}$ that accounts for microscopic heterogeneities. 

Since $\strain^\Omega$ varies throughout the macroscopic domain, a micromodel for $\omega$ is embedded at each Gauss point in $\Omega$ and a microscopic boundary-value equilibrium problem assuming small displacements and strains is solved:
\begin{equation}
\nabla\cdot\stress^\omega=\mbf{0}
\quad\quad
\strain^\omega = \frac{1}{2}\left(\nabla\mbf{u}^\omega+\left(\nabla\mbf{u}^\omega\right)^\T\right)
\end{equation}
\noindent microscopic stress $\stress^\omega$ is related to microscopic strain $\strain^\omega$ with traditional physics-based constitutive models for each phase in the heterogeneous domain. In the general case where the material models feature internal variables $\boldsymbol\alpha$, we can write the constitutive update for the microscale domain as:
\begin{equation}
\model^\omega
\begin{cases}
\state_t^\omega = \mathcal{A}\left(\strain_t^\omega,\state_{t-1}^\omega,\props^\omega\right)\\
\stress_t^\omega = \mathcal{S}\left(\strain_t^\omega,\state_t^\omega,\props^\omega\right)
\end{cases}
\label{EQgeneralmodel}
\end{equation}
\noindent where $\props^\omega$ are the material parameters of the microscopic constituents, the operators $\mathcal{A}$ and $\mathcal{S}$ can be split into an arbitrary number of blocks with different models (\eg elasticity, elastoplasticity, damage) for the different material phases, and $\state^\omega$ is a concatenation of the internal variables of every microscopic Gauss point and therefore fully describes the path-dependent state of the microscopic problem.

In order to determine the strains $\strain^\Omega$ that serve as boundary conditions for the micromodels, a macroscopic small-strain equilibrium problem is solved:
\begin{equation}
\nabla\cdot\boldsymbol{\sigma}^\Omega = \mbf{0}
\quad\quad
\strain^\Omega = \frac{1}{2}\left(\nabla\mbf{u}^\Omega+\left(\nabla\mbf{u}^\Omega\right)^\T\right)
\end{equation}
\noindent but this time no constitutive assumptions are adopted. Macroscale stresses are instead directly homogenized from the microscopic response:
\begin{equation}
\stress^\Omega = \frac{1}{\lvert\omega\rvert}\int_\omega\stress^\omega\mrm{d}\omega
\label{EQfe2stress}
\end{equation}
\noindent which couples the macroscopic strain $\strain^\Omega$ with the microscopic solution. Since \cref{EQfe2disp} also couples the solutions in the opposite direction, a bidirectional coupling is formed which requires the two-scale equilibrium problem to be solved iteratively.

\subsection{Data-driven surrogate modeling}
\label{SEfe2surrogates}

The coupled problem of \cref{SEfe2primer} is extremely computationally demanding. The lower-scale domain $\omega$ usually features complicated geometric features and must therefore be modeled with dense FE meshes in order to ensure accuracy. Worse yet, an independent microscopic problem must be solved at every integration point in $\Omega$ for every iteration of every time step of the simulation. This nested nature quickly forms a computational bottleneck.

Since the bulk of the computational effort lies in solving the micromodels, a popular approach to make multiscale analysis viable for practical applications is to substitute the microscopic FE models by data-driven surrogates. The idea is to perform a number of micromodel simulations under representative boundary conditions and use the resulting stress-strain pairs to train a machine learning model to be deployed when performing the actual two-scale simulations of interest. Naturally, the approach tacitly assumes that the number of offline micromodel computations required to train the model is much smaller than the number of times the microscopic behavior will be computed online. In the following, we use a simple example to demonstrate a number of difficulties associated with training such a model to reproduce path-dependent material behavior.

\subsection{Example: A one-dimensional RNN surrogate}
\label{SEfe2lstm}

\begin{figure}
\centering
\begin{subfigure}[c]{0.45\textwidth}
\centering
\includegraphics[scale=1]{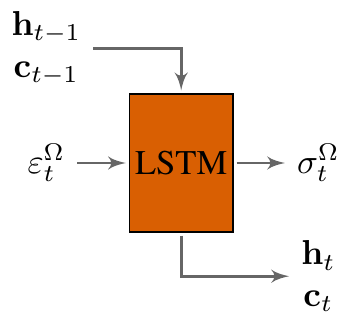}
\caption{Model architecture}
\label{FIlstmscheme}
\end{subfigure}
\begin{subfigure}[c]{0.45\textwidth}
\includegraphics[scale=0.8]{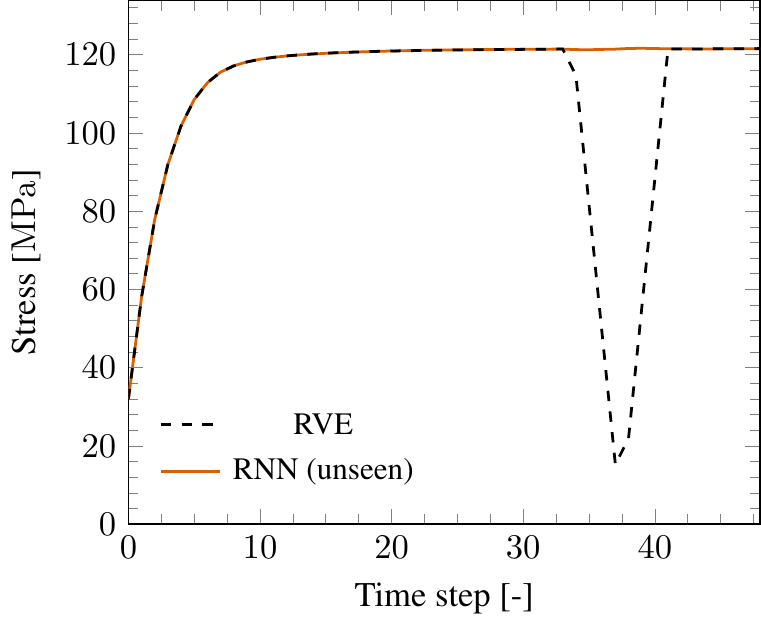}
\caption{Failure to predict unseen unloading}
\label{FIlstmmono}
\end{subfigure}
\caption{An LSTM recurrent neural network as surrogate for 1D path-dependent material behavior trained with only monotonic data.}
\end{figure}

For this demonstration, we train a Long Short-term Memory (LSTM) network \cite{hochreiter_long_1997} to reproduce one-dimensional (single stress/strain component) elastoplasticity. The architecture of the model is shown in \cref{FIlstmscheme} and is implemented in PyTorch \cite{paszke_pytorch_2019}. In order to minimize the risk of overfitting, a pragmatic model selection procedure is performed by first training the model with several non-monotonic strain paths and gradually increasing cell size until reasonable accuracy is obtained. This leads to a parsimonious model with a single LSTM cell with 5 latent units. 

At this point it is interesting to draw a parallel between the network and the micromodel whose behavior is being reproduced: the concatenation of the hidden state $\mbf{h}$ and cell state $\mbf{c}$ of the LSTM cell can be seen as a lower-dimensional surrogate for the set of microscopic internal variables $\state^\omega$ of \cref{EQgeneralmodel}. However, in contrast to the variables in $\state$, the latent variables $\mathbf{h}$ and $\mathbf{c}$ have no physical interpretation and evolve purely according to heuristic memory mechanisms that mimic patterns inferred during training.

First, we train the LSTM using only monotonic data. Since only one strain component is being modeled, this initial dataset is composed simply of one strain path in tension and one in compression. The trained model is then used to predict a tension path with one unloading-reloading cycle. Having never seen unloading during training, the network reverses course and unloads on top of its loading path (\cref{FIlstmmono}). This result is hardly surprising, but sheds light on the potentially deceiving nature of the training procedure: even though we are only concerned with a single strain component, predictions actually take place in an augmented space that describes strain paths in time which can be arbitrarily high-dimensional (as paths can be arbitrarily long).

\begin{figure}
\centering
\begin{subfigure}[b]{0.45\textwidth}
\includegraphics[scale=0.8]{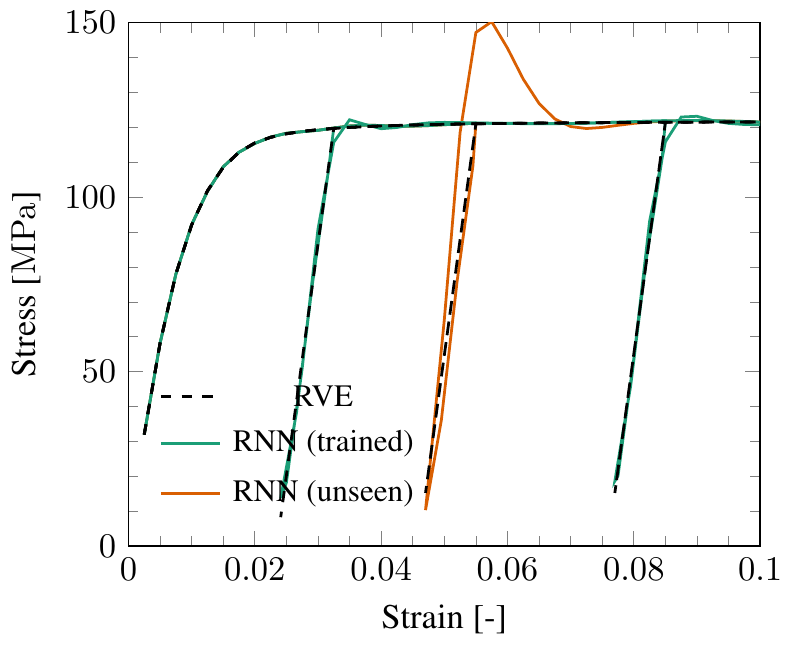}
\caption{Unloading at new strain level}
\label{FIlstmunloading}
\end{subfigure}
\begin{subfigure}[b]{0.45\textwidth}
\includegraphics[scale=0.8]{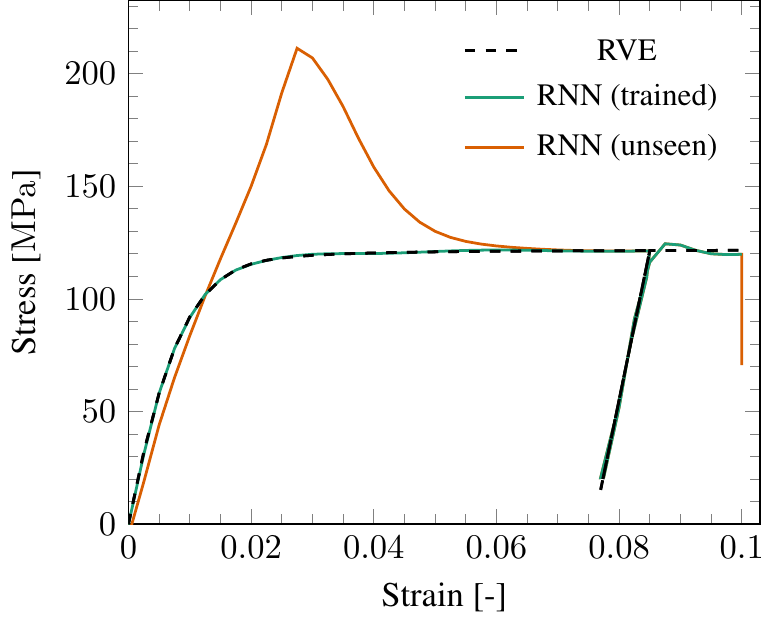}
\caption{Stopping for 5 steps at $\varepsilon=0$ and $\varepsilon=0.1$}
\label{FIlstmstop}
\end{subfigure}
\caption{1D LSTM surrogate trained with unloading/reloading and used to predict unseen unloading paths.}
\label{FIlstmcurse}
\end{figure}

We can further demonstrate this manifestation of the \textit{curse of dimensionality} with the two additional examples of \cref{FIlstmcurse}. In \cref{FIlstmunloading} we train the network with two unloading paths and it fails to predict a third one at an intermediate strain level. Here it can be deceiving to assume the third path can be interpolated from the other two: in the 48-dimensional space of strain paths (we use paths with 48 time steps each) the network is actually operating far away from training data. In \cref{FIlstmstop} the network tries to reproduce a path seen during training but we first let the material rest at zero strain for five time steps before loading starts and for another five time steps at the end of the path. With purely data-driven latent dynamics, the initial rest disturbs the memory structure of the network and causes large deviations for a large portion of the path. For the rest at the end of the path, we see that the surrogate fails to predict the characteristic that the stress does not change upon constant deformation.

Training data-driven models to accurately reproduce path dependency is therefore not straightforward: their latent representations of material state are not interpretable and even phenomena as trivial as resting at zero strain must be learned from data. At the core of successful applications of RNNs to this task are either extensive datasets obtained with carefully crafted sampling strategies \cite{wu_recurrent_2020,logarzoSmartConstitutiveLaws2021} or highly tailored datasets for specific macroscopic problems \cite{ghavamianAcceleratingMultiscaleFinite2019}. Alternatively, active learning frameworks may be used to skip offline training altogether \cite{knap_adaptive_2008,rochaOntheflyConstructionSurrogate2021}, but at the cost of producing slower surrogates.

\section{A hybrid surrogate model}
\label{SEmodel}

In this work we attempt to avoid the curse of dimensionality by relegating to a physics-based material model some of the tasks the RNN of \cref{SEfe2lstm} has to explicitly learn from data. In this section, we further formalize the hybrid approach of \cref{FIapproach} by looking at the roles of each model component and their dependencies in time. We then particularize the model for the case of a feedforward neural network (FNN) encoder and discuss feature selection and numerical stabilization strategies.

\subsection{Evolving material parameters}
\label{SEevolving}

Physics-based material models are traditionally formulated with a fixed set of parameters \props\ either directly computed from a specific set of (numerical) experiments or indirectly from stress-strain measurements in a Maximum Likelihood Estimation (MLE) approach\footnote{The parameters \props\ can also be estimated through Bayesian inference and would therefore be described by a multivariate probability density instead of a fixed set of values. Regardless, that density would still be stationary in time.}. Here we start from the premise that letting (part of) \props\ evolve in time increases flexibility and allows the model to capture more complex material behavior. Conversely, keeping the remainder of the model intact improves interpretability and provides physics-based bias to the data-driven model tasked to learn this evolution.

\begin{figure}
\centering
\includegraphics[scale=1]{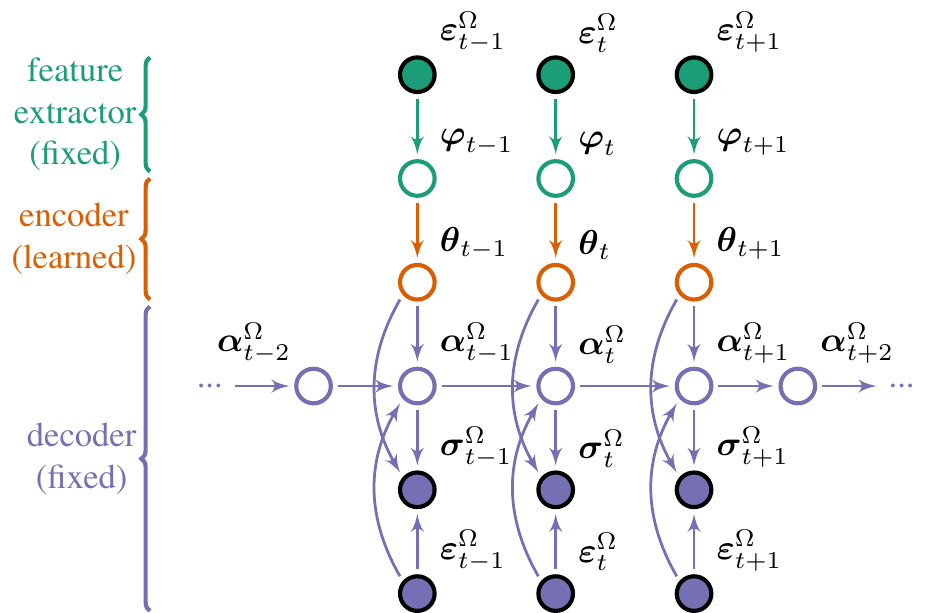}
\caption{Graph representation of the hybrid model architecture combining a data-driven encoder and a physics-based decoder. Filled circles represent observable variables and hollow circles represent latent variables.}
\label{FIgraph}
\end{figure}

In \cref{FIgraph}, the hybrid model of \cref{FIapproach} is unrolled in time for a number of consecutive time steps and represented as a graph showing the dependencies between variables. Filled and hollow nodes represent observed and latent variables, respectively, and are color coded to represent the different model components in \cref{FIapproach}. Similar to the microscale models of \cref{EQgeneralmodel}, we assume the constitutive behavior at the macroscale is given by a physics-based material model:
\begin{equation}
\model^\Omega
\begin{cases}
\state_t^\Omega = \mathcal{A}\left(\strain_t^\Omega,\state_{t-1}^\Omega,\props^\Omega_t\right)\\
\stress_t^\Omega = \mathcal{S}\left(\strain_t^\Omega,\state_t^\Omega,\props^\Omega_t\right)
\end{cases}
\label{EQgeneralmodel2}
\end{equation}
\noindent but now with time-dependent parameters $\props_t$. Note that the model response at time $t$ depends on the material state at time $t-1$ through a set of internal variables $\state^\Omega_{t-1}$ (\cref{FIgraph}). This gives the model a recurrent nature not unlike that of the RNN of \cref{FIlstmscheme} with its state variables $\mathbf{c}$ and $\mathbf{h}$. The advantage here is that \state\ has clear physical interpretation (plastic strains, damage variables, etc) and its evolution is handled by the fixed operator $\mathcal{A}$ composed of clearly interpretable algorithmic steps grounded in physics and/or classical material phenomenology (\eg a return mapping algorithm).

On the encoder side, we let the material properties \props\ evolve according to an evolution operator \data\ whose shape is learned from data:
\begin{equation}
\props_t= \data\left(\boldsymbol\varphi_t\right)
\end{equation}
\noindent as a function of a set of features \features\ that are themselves obtained from the macroscopic strains through a feature extractor \extractor:
\begin{equation}
\boldsymbol\varphi_t = \extractor\left(\boldsymbol\varepsilon^\Omega_t\right)
\end{equation}
\noindent where $\features_t$ could be simply the strains themselves or other quantities derived from it. More importantly, note that $\props_t$ depends only on the current features $\features_t$ and we therefore assume the encoder is not recurrent (\cref{FIgraph}). This choice effectively limits the flexibility of \data\ and makes the hybrid surrogate fully rely on the more robust model $\model^\Omega$ to explain path-dependent phenomena, helping counter the curse of dimensionality associated with sampling strain paths. For instance, it opens up the possibility to train the surrogate exclusively with monotonic data, as we will demonstrate in the examples of \cref{SEexamples}.

In the following sections, we particularize the model for the case of \data\ being a fully-connected neural network and for specific choices of \extractor\ and \model. Nevertheless, the general architecture of \cref{FIapproach,FIgraph} is meant to be as flexible as possible:
\begin{itemize}
\item The nature and dimensionality of \features\ is not tied to that of $\strain^\Omega$ since strains are also given directly to $\model^\Omega$;
\item Other machine learning models for regression can also be used as \data, and it could in principle be split into different models handling the evolution of different subsets of \props. Any number of model parameters may also be left out of \props\ and either fixed as constants or optimized to constant values during training;
\item No assumption is made on the form of $\model^\Omega$ or the nature or dimensionality of $\state^\Omega$. Instead of a single model, it could also for instance be a mixture of physics-based models combined with analytical homogenization techniques.
\end{itemize}

\subsection{Feature extractors}
\label{SEfeatureselection}

A pragmatic choice for \extractor\ is to simply assume \features\ is the macroscopic strain vector $\strain^\Omega$ itself. It is also a familiar one, as we can then relate the resulting model to conventional surrogates mapping strains to stresses. However, since macroscopic strains are also directly passed on to the decoder, the architecture gives us the freedom to experiment with different features.

\begin{figure}
\centering
\begin{subfigure}[b]{0.45\textwidth}
\centering
\includegraphics[scale=1]{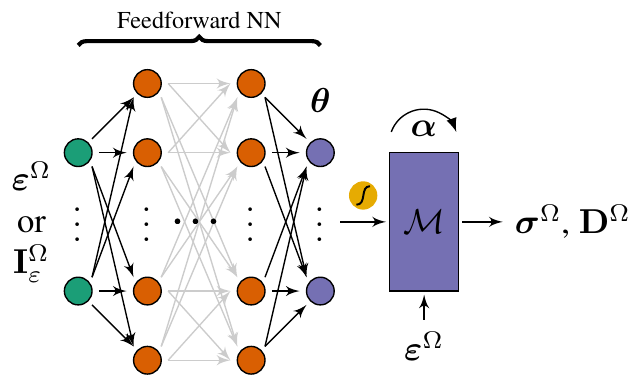}
\caption{Macroscopic strains or invariants directly used as features}
\label{FIarch1}
\end{subfigure}
\begin{subfigure}[b]{0.45\textwidth}
\centering
\includegraphics[scale=1]{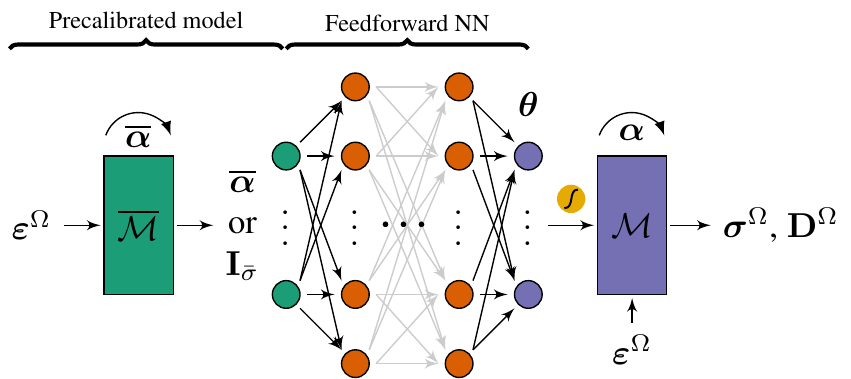}
\caption{Features extracted from an informative microscopic constitutive model}
\label{FIarch2}
\end{subfigure}
\caption{The two types of FNN-based model architectures explored in this work, with different feature extraction steps.}
\label{FIarchitectures}
\end{figure}

\cref{FIarchitectures} shows the two model architectures we explore in this work. For the two variants in \cref{FIarch1} we either use $\strain^\Omega$ itself or a set of small-strain invariants of the macroscopic strain tensor of increasing dimensionality:
\begin{equation}
\mbf{I}^\Omega_\varepsilon = \begin{bmatrix} I^\varepsilon_1\end{bmatrix}
\quad\mrm{or}\quad
\mbf{I}^\Omega_\varepsilon = \begin{bmatrix} I^\varepsilon_1 & I^\varepsilon_2\end{bmatrix}
\end{equation}
\noindent where the variants are given by the well-known expressions:
\begin{equation}
I^\varepsilon_1 = \tr\left(\boldsymbol\varepsilon\right), \quad
I^\varepsilon_2 = \frac{1}{2}\left(\tr\left(\boldsymbol\varepsilon\right)^2-\tr\left(\boldsymbol\varepsilon^2\right)\right)
\end{equation}
\noindent Additionally, since the current study focus on elastoplasticity, it is also interesting to explore feature spaces including invariants from the deviatoric strain tensor:
\begin{equation}
\mbf{I}^\Omega_\varepsilon = \begin{bmatrix} J^\varepsilon_2 \end{bmatrix}
\quad\mrm{or}\quad
\mbf{I}^\Omega_\varepsilon = \begin{bmatrix} I^\varepsilon_1 & J^\varepsilon_2 \end{bmatrix}
\end{equation}
\noindent where:
\begin{equation}
J^\varepsilon_2 = \frac{1}{3}\left(I^{\varepsilon}_1\right)^2 - I^\varepsilon_2
\end{equation}
\noindent By using features based on invariants and since the decoder material model is itself already frame invariant for small strains, it follows that the resulting surrogate will naturally inherit this beneficial characteristic. This stands in contrast with traditional black-box surrogates mapping strains to stresses. Furthermore, opting for invariant-based features can be seen as a physics-based dimensionality reduction operation that can potentially reduce the amount of data needed to train the hybrid model.

We also investigate the possibility of extracting features from the outputs of a precalibrated physics-based material model $\overline\model$ subjected to the same strain path seen at the macroscale (\cref{FIarch2}). Note that this specific architecture introduces an additional recurrent component to the model through the set $\overline\state$ of internal variables of $\overline\model$. From a machine learning perspective, the role of $\overline\model$ would be analogous to that of a temporal convolution operator or an RNN cell appended to the encoder. The key difference, however, is that $\overline\model$ is fixed \textit{a priori} and therefore should not require extra sampling effort with respect to the more straightforward extractor in \cref{FIarch1}.

Naturally, different choices for $\overline\model$ yield models with distinct learning capabilities, and we therefore assume $\overline\model$ encapsulates relevant information about not only the current values of $\strain^\Omega$ but also of its history. In the present scenario where the data is coming from micromodel computations, we opt for the intuitive choice of having $\overline\model$ be one of the known constitutive models used to describe the microscopic material phases. We can therefore conceptually see $\overline\model$ as an imaginary representative material point at the microscale that is always subjected to the average micromodel strain. We then use either a subset of its internal variables $\overline\state$ or a set of invariants $\mbf{I}_{\bar{\boldsymbol{\sigma}}}$ of its stress outputs as features.

\subsection{Neural network encoder}
\label{SEmlps}

For simplicity, we opt for modeling the evolution of \props\ using classical feedforward neural networks with fully-connected layers. As both architectures in \cref{FIarchitectures} ultimately compute macroscopic stresses given macroscopic strains, we can use supervised learning to train the model with a straightforward Maximum Likelihood approach. Gathering the complete set of network weights in a vector $\weights$ and seeing the complete surrogate as a monolithic model that computes an approximation $\widehat{\stress}$ for stresses, we adopt the following observation model for the snapshot stresses \stress: 
\begin{equation}
\stress = \widehat{\stress}\left(\strain,\weights\right) + \xi, \quad \xi\sim\mathcal{N}\left(\xi\vert\mbf{0},\beta^{-1}\mbf{I}\right)
\end{equation}
\noindent where the superscript $\Omega$ is dropped for convenience, $\mbf{I}$ is an identity matrix, and $\xi$ is an additive Gaussian noise\footnote{Even though our observations come from a computer model and can be considered noiseless, the surrogate $\widehat\stress$ is in general not arbitrarily flexible and the random variable $\xi$ is therefore still necessary to explain why the model does not exactly fit every single observation in the dataset.}. Under the assumption of a squared loss, maximizing the likelihood of a training dataset with $N$ observations amounts to minimizing the loss function \cite{bishopPatternRecognitionMachine2006}:
\begin{equation}
L = \frac{1}{2}\sum_{n=1}^N\norm{\stress_n - \widehat{\stress}\left(\strain_n,\weights\right)}^2
\label{EQloss}
\end{equation}
\noindent with the variance of the noise that explains data misfit being simply $\beta=N/2L$. The resulting loss function is the same one used for conventional data-driven surrogates and is therefore straightforward to implement. 

Nevertheless, it is worth noting that since we cannot directly observe \props, computing the gradients of $L$ with respect to $\weights$ involves backpropagating derivatives through the decoder \model. Furthermore, since $\weights$ affects the evolution of the internal variables $\state$, backpropagation in time becomes necessary. Starting from \cref{EQloss} and walking back through the graph of \cref{FIgraph}, the gradient of the loss at time step $t$ of a given strain path is given by:
\begin{equation}
\frac{\partial L_t}{\partial\weights} =\frac{\partial L}{\partial\widehat{\stress}_t}\left\{
\frac{\partial\widehat{\stress}_t}{\partial\props_t}\frac{\partial\props_t}{\partial\weights}
+ \frac{\partial\widehat{\stress}_t}{\partial\state_t}\frac{\partial\state_t}{\partial\props_t}\frac{\partial\props_t}{\partial\weights}
+ \frac{\partial\widehat{\stress}_t}{\partial\state_t}\sum_{\bar{t}=t-1}^{1}\left[\left(\prod_{\tilde{t}=t}^{\bar{t}+1}\frac{\partial\state_{\tilde{t}}}{\partial\state_{\tilde{t}-1}}\right)\frac{\partial\state_{\bar{t}}}{\partial\props_{\bar{t}}}\frac{\partial\props_{\bar{t}}}{\partial\weights}\right]
\right\}
\end{equation}
\noindent where the remaining gradient chain $\partial\props/\partial\weights$ is computed with conventional backpropagation through the network. If \model\ is implemented in a code base that allows for automatic differentiation (\eg in PyTorch), these time dependencies are naturally taken into account as long as a persistent gradient tape is used within each strain path\footnote{This is already the case for RNNs, so switching from RNNs to the present model should require little to no changes to the way training is performed.}. In this work we instead implement network training directly into an existing FE code, and therefore opt for the pragmatic approach of computing all partial derivatives of quantities derived from \model\ using finite differences.

Finally, in order to enforce upper and lower bounds for \props\ and avoid unphysical parameter values (\eg negative elasticity moduli), we apply sigmoid activation to the final layer of the network and scale the parameters back from a $[0,1]$ range using predefined bounds:
\begin{equation}
\theta_i = \theta_i^\mrm{low} + \theta_i^\sigma\left(\theta_i^\mrm{upp}-\theta_i^\mrm{low}\right)
\label{EQbounds}
\end{equation}

\subsection{Material decoders}
\label{SEdecoders}

As previously mentioned, any constitutive model can in principle be used as \model. For the present study we focus on reproducing elastoplasticity and therefore narrow our choices down to the following set of potential decoders with increasing levels of complexity. The simplest one is a linear-elastic isotropic material with no internal variables: 
\begin{equation}
\sigma_{ij} = D_{ijkl}\varepsilon_{kl} \quad\mrm{with}\quad D_{ijkl} = G\left(\delta_{ij}\delta_{kl} + \delta_{il}\delta_{jk}\right) + \left(K-\frac{2}{3}G\right)\delta_{ij}\delta_{kl}
\label{EQlineardecoder}
\end{equation}
\noindent where index notation is used for convenience. For this model, \props\ comprises only the bulk and shear moduli $K$ and $G$, or equivalently the Young's modulus $E$ the Poisson's ratio $\nu$.

The second decoder option is a simple plasticity model with $J_2$ (von Mises) flow. The stress update in this case becomes:
\begin{equation}
\sigma_{ij} = D_{ijkl}\left(\varepsilon_{ij} - \varepsilon^\p_{ij}\right)
\label{EQplasticstressupdate}
\end{equation}
\noindent where strain is additively decomposed into elastic and plastic ($\boldsymbol{\varepsilon}^\p$) contributions. The yield criterium and plastic flow rule are given by:
\begin{equation}
\phi=\sqrt{3J_2^\sigma} - \sigma_\mrm{y} \leq 0
\quad\mrm{and}\quad
\Delta\varepsilon^\p_{ij} = \Delta\gamma\sqrt{\frac{3}{2}}\frac{S_{ij}}{\norm{S_{ij}}_\mathrm{F}}
\label{EQj2}
\end{equation}
\noindent where $\mbf{S}$ is the deviatoric part of the stresses, $\gamma$ is a plastic multiplier, $\sigma_\mrm{y}$ is a yield stress parameter and we write the Frobenius norm as $\norm{\cdot}_\mrm{F}$. In order to keep the model as simple as possible, we assume $\sigma_\mrm{y}$ is a material constant and therefore end up with a perfectly-plastic model with associative flow. The internal variables of this model are components of the plastic strain vector $\boldsymbol{\varepsilon}^\p$ and the only new material parameter is the yield stress $\sigma_\mrm{y}$.

Finally, we also consider the more complex pressure-dependent, non-associative plasticity model proposed by Melro \etal \cite{melro_micromechanical_2013}. Stress update is the same as in \cref{EQplasticstressupdate}, but yield surface and plastic flow are given by:
\begin{equation}
\phi=6J^\sigma_2 + 2I^\sigma_1\left(\sigma_\comp-\sigma_\tens\right) - 2\sigma_\comp\sigma_\tens\leq 0
\quad\mrm{and}\quad
\Delta\varepsilon^\p_{ij} = \Delta\gamma\left(3S_{ij}+\frac{1-2\nu_\p}{1+\nu_\p}I_1^\sigma\delta_{ij}\right)
\label{EQmelro}
\end{equation}
\noindent where $\delta_{ij}$ is the Kronecker delta, $\sigma_\tens$ and $\sigma_\comp$ are yield stresses in tension and compression, respectively, and $\nu_\p$ is a new parameter controlling plastic contraction and allowing for compressible plastic flow. Hardening can be described by making the yield stresses general functions of $\strain^\p$, but when used as a decoder we assume $\sigma_\tens$ and $\sigma_\comp$ do not depend on $\strain^\p$ and instead let the decoder \data\ describe their evolution.

The model by Melro \etal \cite{melro_micromechanical_2013} is also the one used to describe the microscopic material phase responsible for the nonlinear behavior observed when homogenizing micromodel response, and can therefore be seen as the natural choice for \model. Nevertheless, the other two decoders can provide interesting insights on the effect of introducing different levels of bias to the hybrid model.

\subsection{Online predictions and inherited stability}
\label{SEstabilization}

The architecture of \cref{FIapproach} is developed to be minimally intrusive and allow for existing material models to be used as decoders with minimum effort. We therefore implement the online routine of the model as a wrapper around an existing implementation of \model. The basic structure of the wrapper can be seen in \cref{ALwrapper}.
The hybrid nature of the model allows for a robust approach that ensures the numerical stability of the original model \model\ is inherited by the surrogate. This is achieved by only updating \props\ at the end of each time step, after the global implicit Newton-Raphson scheme converges. Material properties are therefore fixed while the global solver is iterating, and that means the tangent stiffness $\mbf{D}$ comes directly from \model\ and inherits its stability features. 

\begin{algorithm}
\caption{Material wrapper implementing the online component of the hybrid surrogate.}
\label{ALwrapper}
\BlankLine
\KwIn{strain $\strain_\mrm{new}^\Omega$ at macroscopic Gauss point}
\KwOut{stress $\stress^\Omega$ and stiffness $\mbf{D}^\Omega$ at macroscopic Gauss point}
\BlankLine
use nested model with converged parameters and internal state: $\left(\stress^\Omega,\mbf{D}^\Omega,\state_\mrm{new}\right) \leftarrow \model\left(\strain_\mrm{new}^\Omega,\state_\mrm{old},\props\right)$\;
\If{global solver has converged}{
store latest converged strain: $\strain_\mrm{old}\leftarrow\strain_\mrm{new}$\;
commit material history: $\state_\mrm{old}\leftarrow\state_\mrm{new}$\;
compute new features: $\features_\mrm{new}\leftarrow\extractor\left(\strain_\mrm{new}\right)$\;
update model parameters for the upcoming time step: $\props\leftarrow\data\left(\features_\mrm{new}\right)$\;
}
\If{first global iteration of time step \And Gauss point is unstable}{
stabilize encoder: $\data\leftarrow\mathtt{stabilizeNetwork}\left(\strain_\mrm{new}^\Omega\right)$\;
recompute features: $\features_\mrm{old}\leftarrow\extractor\left(\strain_\mrm{old}^\Omega\right)$\;
recompute model parameters for the current time step: $\props\leftarrow\data\left(\features_\mrm{old}\right)$\;
}
\Return{$\bs\sigma^\Omega,\mbf{D}^\Omega$}
\end{algorithm}


As an example, the $J_2$ plasticity model of \cref{EQj2} is unconditionally stable as long as its hardening modulus $h\geq 0$ for any $\left(\strain^\Omega_t,\state^\Omega_t\right)$, which is the case for the perfectly-plastic version we consider here. It then follows that any hybrid surrogate with $J_2$ decoder is also unconditionally stable. Note that this is only possible because strains are directly passed on to the decoder and would therefore not be an option for conventional surrogates (\eg the RNN of \cref{FIlstmcurse}). For those surrogates, the tangent stiffness would come directly from the jacobian of a highly-flexible data-driven model, often at the cost of numerical stability.

\subsection{Numerical stabilization}
Nevertheless, the decoder \model\ may be inherently unstable even with fixed material constants. This is for instance the case for the model by Melro \etal \cite{melro_micromechanical_2013}: the non-associative flow rule of \cref{EQmelro} can cause the tangent stiffness $\mbf{D}^\Omega$ to lose positive definiteness under certain strain conditions and for certain combinations of model parameters. To accommodate such a scenario and open up the possibility for online model adaptivity in other contexts, we propose a scheme for updating the encoder \data\ on the fly in order to enforce extra constraints locally. 

Back to \cref{ALwrapper}, at the beginning of a new time step we keep \props\ fixed to the one obtained with converged strains from the previous step and let the solver make a first strain prediction. After this first iteration, a stability criterion is checked and used to define a new loss function that can be used to update network weights in case instability is detected. Here we employ the determinant of the acoustic tensor $\mbf{Q}$:
\begin{equation}
  \mbf{Q} = \mbf{n}_d^\T\mbf{D}^\Omega\mbf{n}_d
\end{equation}
\noindent where $\mbf{n}_d$ is the vector normal to the strain localization direction creating the instability, which we find through an angle sweep procedure as in \cite{van_der_meer_computational_2022}. We then use $\det\left(\mbf{Q}\right)$ as a metric of stability and trigger a retraining procedure in case a negative value is detected. We then introduce a new loss function:
\begin{equation}
L_\mrm{Q} = -\frac{\left<\det\left(\mbf{Q}\right)\right>_-}{\det\left(\mbf{Q}_0\right)}
\label{EQstabilityloss}
\end{equation}
\noindent where $\left<\cdot\right>_{-}$ extracts the negative part of its operand and $\mbf{Q}_0$ is a reference value for the acoustic tensor computed at the start of the simulation. We minimize this new loss at every unstable point for a small number of epochs with low learning rate, and to discourage significant drifts from the original model we finish the stabilization procedure by updating the network using the original loss of \cref{EQloss} for a single minibatch. Finally, \props\ is updated using the retrained model and is kept fixed for the remaining iterations\footnote{Changing \data\ and therefore \props\ after every iteration would not work in favor of improving stability, but rather have the opposite effect.}. Note that the local constraint of \cref{EQstabilityloss} is therefore only enforced in a soft way and remaining instabilities might still cause the global solver to diverge, in which case we cancel the current increment, go back to the beginning of the time step and allow for the procedure to be triggered again.

\section{Numerical examples}
\label{SEexamples}

The proposed model was implemented in an in-house Finite Element code developed using the open-source C++ numerical analysis library Jem/Jive \cite{nguyen-thanhJiveOpenSource2020}. In order to allow for seamless online retraining, network training was also implemented within the same code. We start this section by describing the datasets and model selection strategies used to build the surrogates. We then investigate the performance of the approach under several choices of encoders and decoders. Finally, we use the model within an \fetwo\ simulation and demonstrate the online stabilization procedure of \cref{SEstabilization}. All simulations are performed on cluster nodes equipped with Xeon E5-2630V4 processors and \SI{128}{GB} RAM running CentOS 7.

\subsection{Data sampling and model selection}

\begin{figure}
  \centering
  \includegraphics[scale=1]{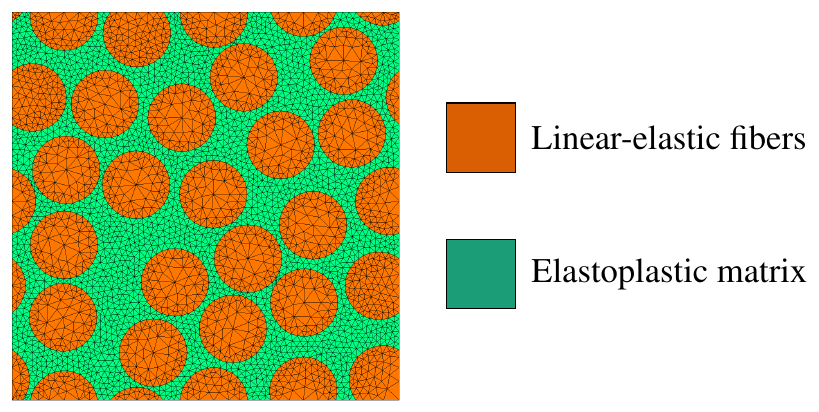}
  \caption{The micromodel used in the examples of this work.}
  \label{FImesh}
\end{figure}

Models are trained to reproduce the behavior of the fiber-reinforced composite micromodel shown in \cref{FImesh}. Fibers are modeled as linear-elastic and the matrix is described by the pressure-dependent non-associative elastoplastic model by Melro \etal \cite{melro_micromechanical_2013} (\cref{SEdecoders}). Microscale material properties are adopted from \cite{vandermeerMicromechanicalValidationMesomodel2016}. The microscopic geometry shown in \cref{FImesh} results from an RVE study performed in \cite{vandermeerMicromechanicalValidationMesomodel2016} and is therefore considered representative. Following the discussion in \cref{SEmodel}, our aim is to investigate up to which extent it is possible to circumvent the curse of dimensionality associated with path dependency by training surrogates exclusively on monotonic strain paths and having time-dependent behavior arise naturally from a physics-based decoder. We therefore limit ourselves to monotonic paths for training. For consistency, we also employ exclusively motononic data to perform model selection.

\begin{figure}
  \centering
  \begin{subfigure}[c]{0.33\textwidth}
    \includegraphics[scale=0.6]{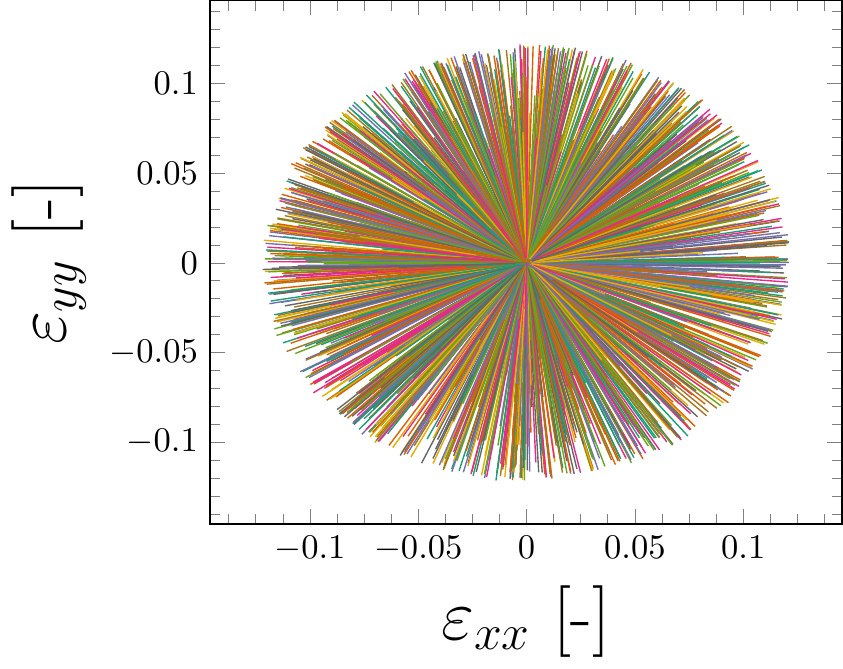}
  \end{subfigure}
  \begin{subfigure}[c]{0.33\textwidth}
    \includegraphics[scale=0.6]{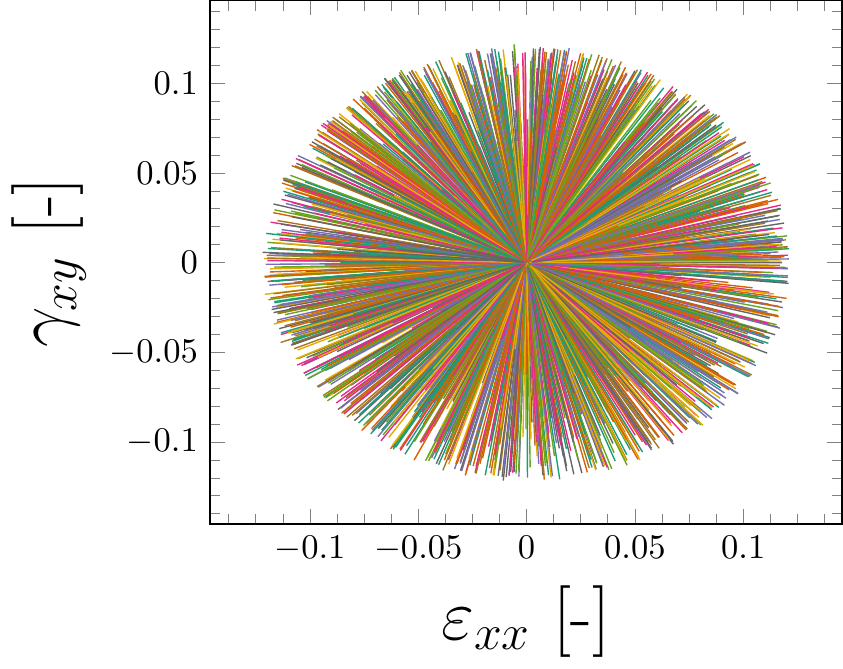}
  \end{subfigure}
  \begin{subfigure}[c]{0.33\textwidth}
    \includegraphics[scale=0.6]{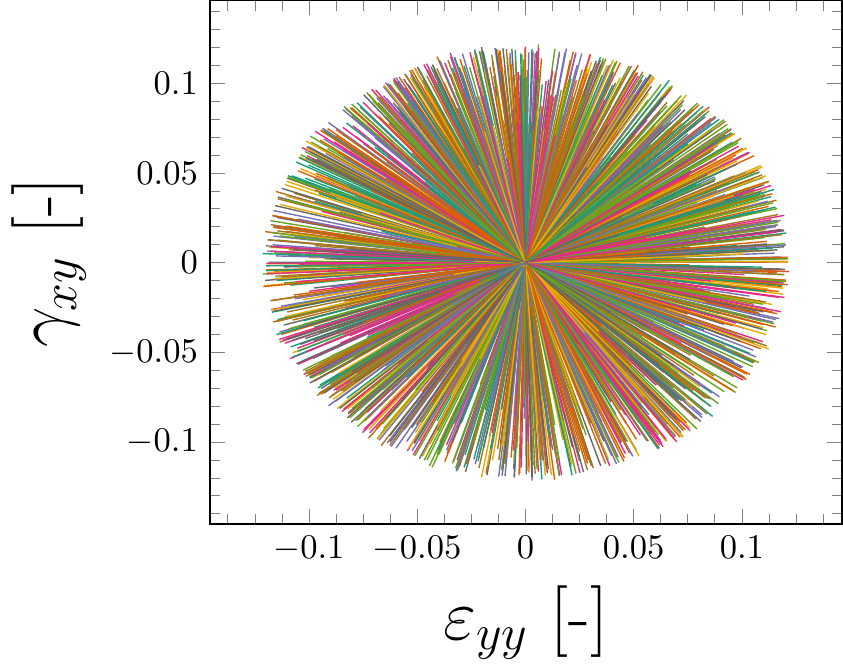}
  \end{subfigure}
  \\
  \begin{subfigure}[c]{0.33\textwidth}
    \includegraphics[scale=0.6]{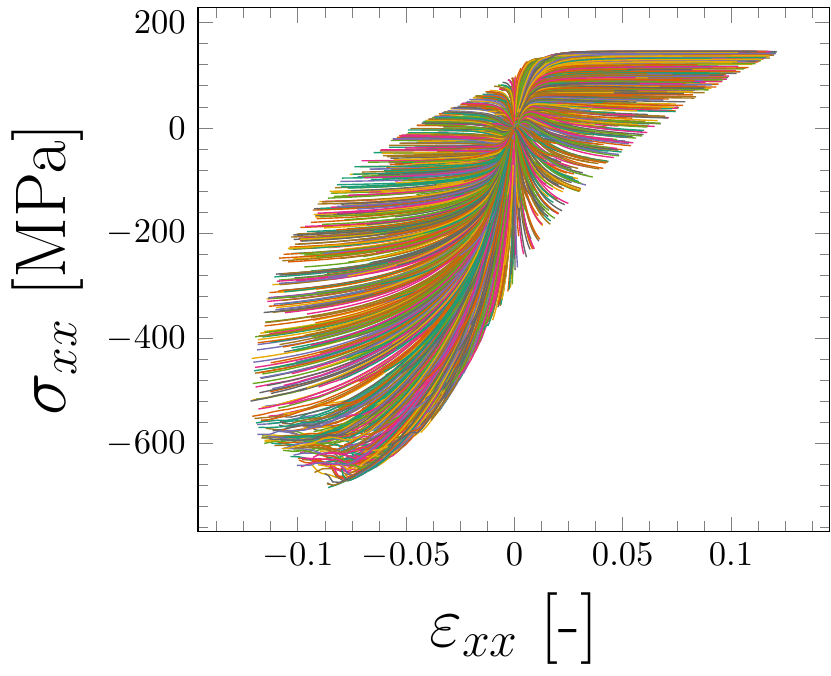}
  \end{subfigure}
  \begin{subfigure}[c]{0.33\textwidth}
    \includegraphics[scale=0.6]{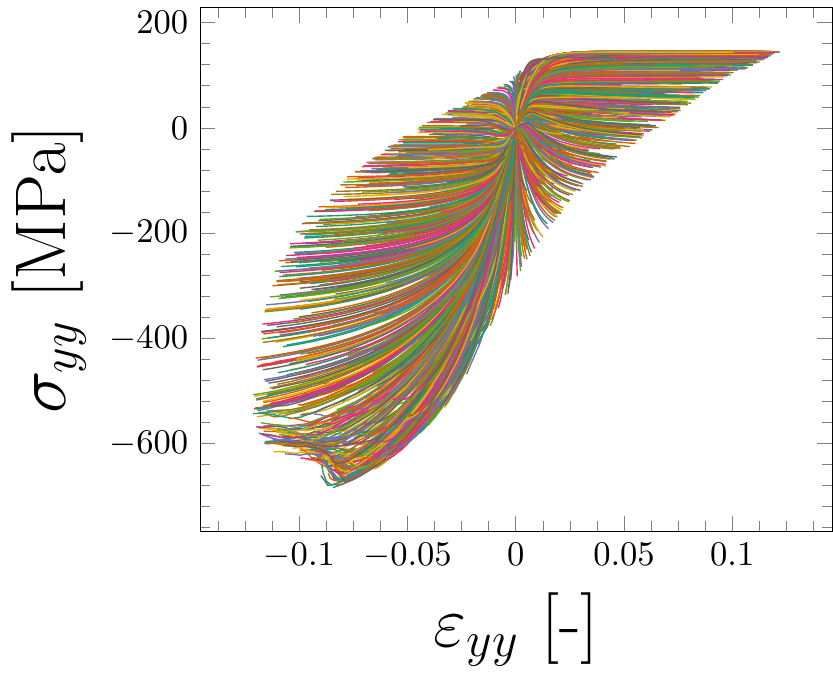}
  \end{subfigure}
  \begin{subfigure}[c]{0.33\textwidth}
    \includegraphics[scale=0.6]{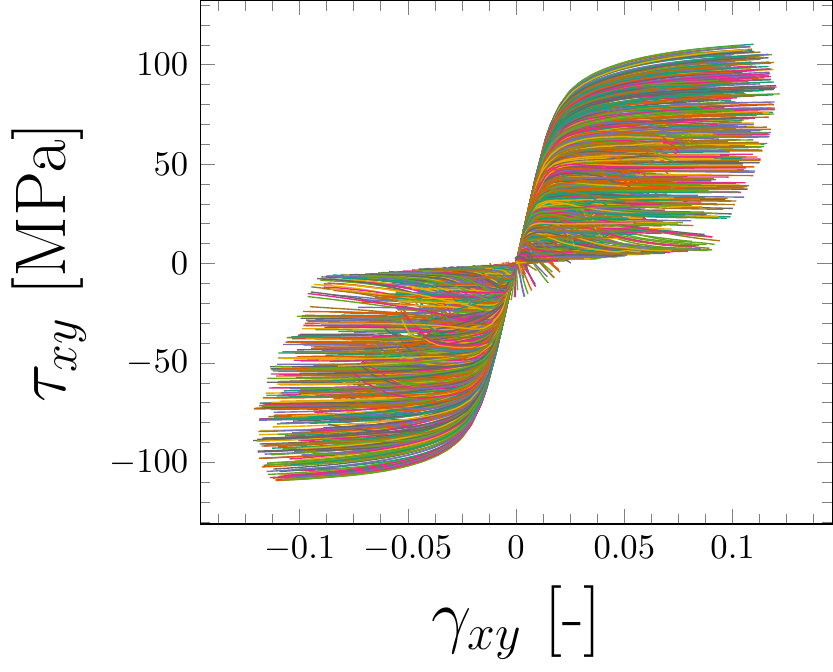}
  \end{subfigure}
  \caption{The complete plane strain dataset used to train the surrogates, comprising \SI{2000}{} monotonic strain-stress paths. A similar dataset is generated under plane stress conditions.}
  \label{FIdataset}
\end{figure}

For efficiency, we limit the present investigation to 2D simulations (\ie three strain components) in the plane perpendicular to the fibers, but nevertheless expect the discussion and conclusions to generalize to 3D simulations as long as appropriate orthotropic decoders are employed. Datasets with \SI{2000}{} monotonic strain paths are generated under both plane strain and plane stress assumptions. \cref{FIdataset} shows the complete plane strain dataset, with a similar one also being generated for plane stress. Each path is generated with an \fetwo\ simulation of a single macroscopic element under displacement control along a fixed direction in strain space sampled from a uniform distribution. To circumvent convergence issues, we employ an adaptive time stepping technique that progressively reduces time step size when the simulation does not converge and gradually increases it back for subsequent increments. The simulations are stopped once a strain norm of \SI{10}{\percent} is reached. As the adaptive scheme leads to paths with different numbers of time increments, we balance the dataset by ensuring every path is composed of \SI{30}{} steps with strain norms as equally spaced as possible.

To keep model selection straightforward and avoid the need for cumbersome k-fold cross validation or bootstrapping, we train a preliminary model with enough flexibility and an extensive training dataset and gradually increase the size of the validation set until the validation error converges to a good estimate of the expected prediction error \cite{hastie}. This results in validation sets with \SI{500}{} paths selected at random from the original datasets, leaving \SI{1500}{} paths to be used for training. We then perform model selection by gradually increasing the complexity of our FNN encoders until the validation error stabilizes. From experimenting with different architectures, we find that encoders with 5 hidden layers of 50 units each with Scaled Exponential Linear Unit (SELU) \cite{klambauerSelfNormalizingNeuralNetworks2017} activation provide enough flexibility for all the examples treated here. To ensure enough regularization when computing learning curves with small datasets, we employ Bernoulli dropout layers with a rate of \SI{1}{\percent} after every hidden layer. Networks are trained for \SI{20000}{} epochs and the model with lowest historical validation error is kept after training, further reducing the risk of overfitting on small datasets.

\begin{figure}
  \centering
  \begin{subfigure}[c]{0.33\textwidth}
    \includegraphics[scale=0.65]{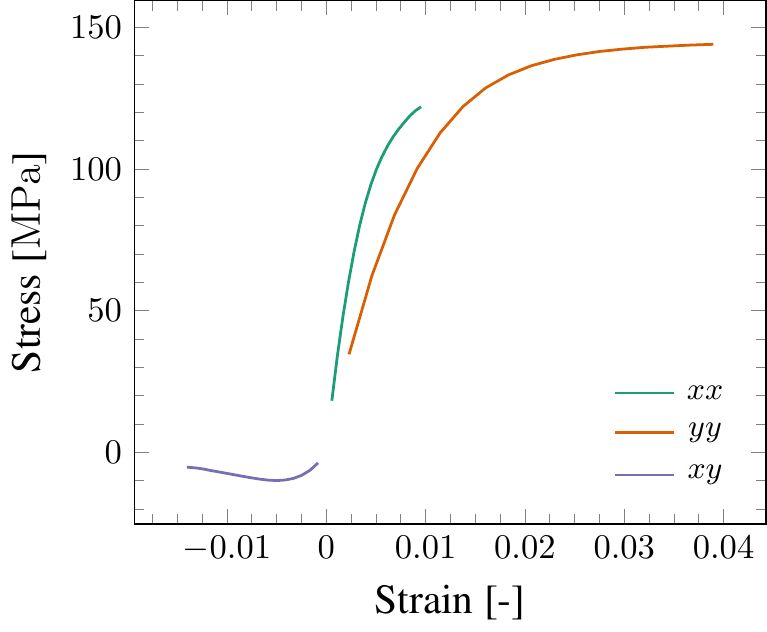}
    \caption{Monotonic}
  \end{subfigure}
  \begin{subfigure}[c]{0.33\textwidth}
    \includegraphics[scale=0.65]{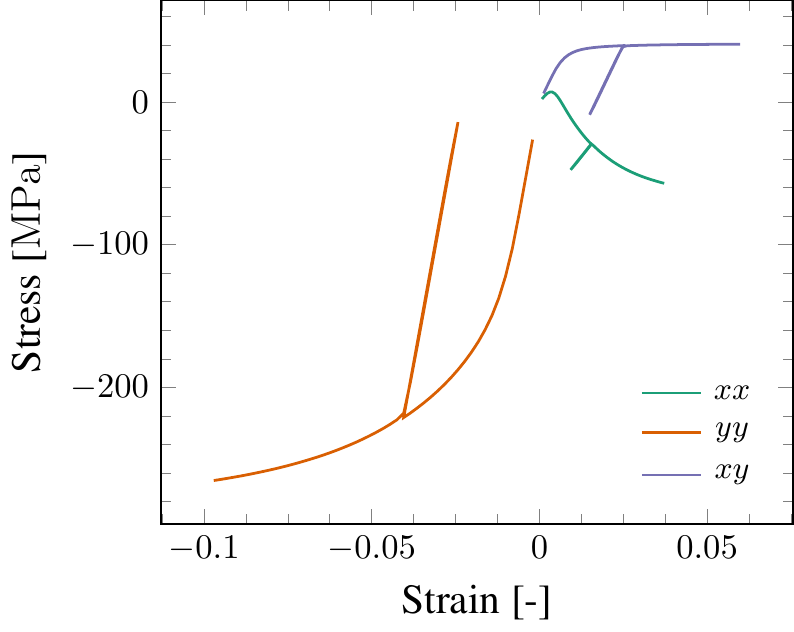}
    \caption{Unloading-reloading}
  \end{subfigure}
  \begin{subfigure}[c]{0.33\textwidth}
    \includegraphics[scale=0.65]{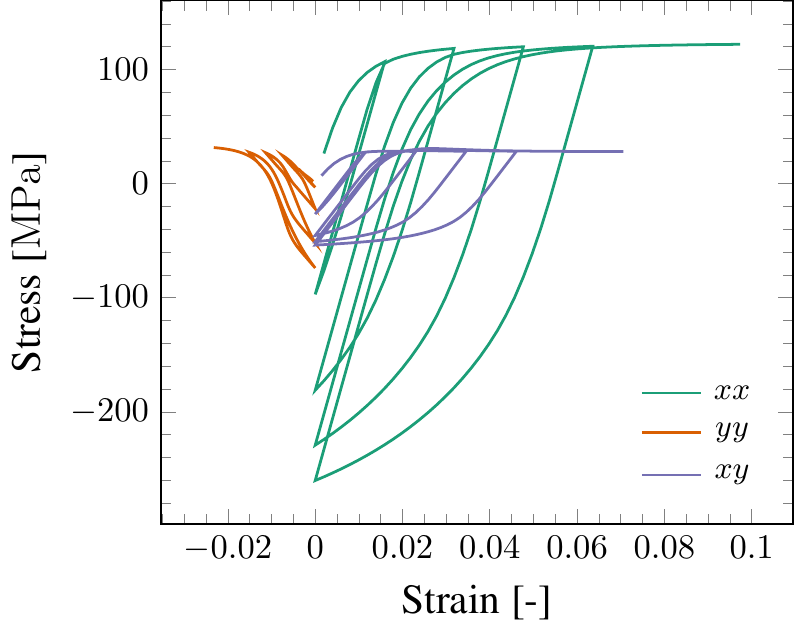}
    \caption{Slow cycling}
    \label{FItestexamplescycles}
  \end{subfigure}
  \caption{Examples from a test dataset with 50 paths of each type. They are not used to train any of the networks or perform model selection.}
  \label{FItestexamples}
\end{figure}

To assess the capabilities of the trained surrogates, we compute an additional test dataset comprising \SI{50}{} monotonic, \SI{50}{} unloading-reloading and \SI{50}{} slow cycling paths, examples of which are shown in \cref{FItestexamples}. To keep the comparisons fair, none of these paths are used to perform model selection and are therefore only considered after the surrogates are trained. We will use example curves like those from \cref{FItestexamples} for visual inspection of the model performance, but also the complete sets of \SI{50}{} curves each for more rigorous statistical analysis.

\subsection{Elastic decoder}

It is interesting to first consider the simple linear-elastic decoder of \cref{EQlineardecoder}, as it has no internal variables and therefore leads to a surrogate model comparable in nature to a conventional FNN trained on stress-strain pairs. As we will demonstrate, however, the limited physical bias provided by such simple model already proves advantageous. Here we let both elastic properties be controlled by the learned encoder:
\begin{equation}
  \props = \begin{bmatrix}E&\nu\end{bmatrix}
\end{equation}
\noindent where the bounds $10^1<E<10^5$ and $0<\nu<0.5$ are enforced as described in \cref{EQbounds}.

\begin{figure}
  \centering
  \begin{subfigure}[c]{0.4\textwidth}
    \includegraphics[scale=0.8]{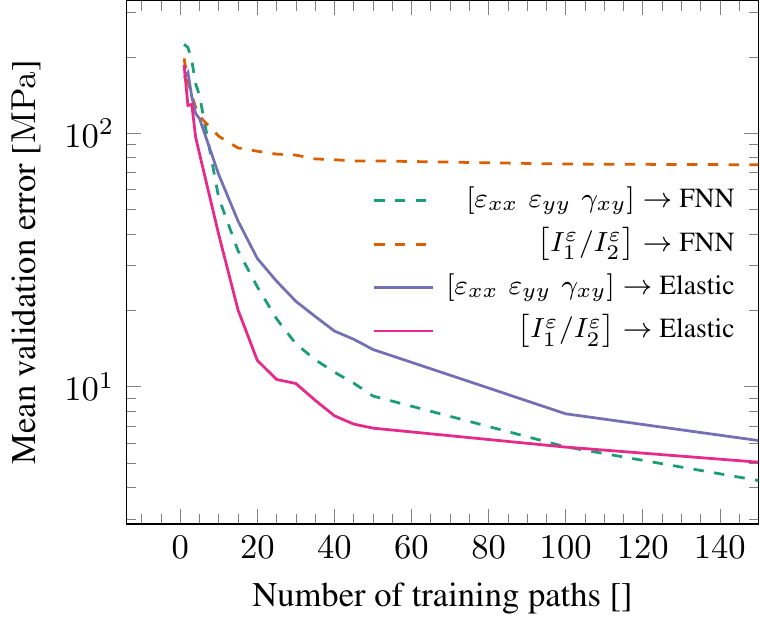}
    \caption{Expectations over 50 datasets of each size}
    \label{FIelasticvsfnn1}
  \end{subfigure}
  \begin{subfigure}[c]{0.4\textwidth}
    \includegraphics[scale=0.8]{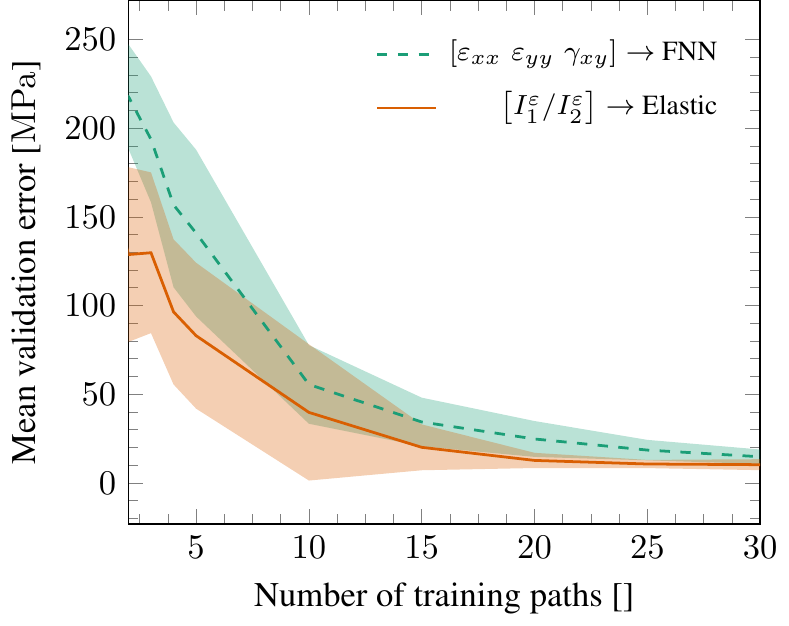}
    \caption{Detailed comparison including standard deviations}
    \label{FIelasticvsfnn2}
  \end{subfigure}
  \caption{Learning curves of models with elastic decoders and conventional FNN models. Mean error over the 500 validation monotonic paths.}
  \label{FIelasticvsfnn}
\end{figure}

We first perform a feature selection study and investigate how efficiently the model learns as the size of the dataset is increased. From the original plane strain training dataset of \SI{1500}{} monotonic strain paths, we draw datasets with sizes ranging between \SI{1}{} and \SI{150}{} paths without replacement and use them to train networks with different encoder features. To get a reliable estimate of the expected prediction error, we repeat this process \SI{50}{} times for each dataset size and encoder type, and for comparison we also do the same for conventional FNNs trained directly on stress targets (keeping the same architecture but going directly from the final hidden layer to stresses). This amounts to a total of \SI{3400}{} trained networks from which we can compute an estimate of the prediction error by averaging $\norm{\boldsymbol{\sigma}-\widehat{\boldsymbol{\sigma}}}$ over the \SI{500}{} paths left for validation. 

\cref{FIelasticvsfnn1} plots averages of the validation error over the \SI{50}{} training datasets used for each size. Although the hybrid architecture does not show an advantage over the FNN when the encoder is trained on strain features, there is a clear gain in learning speed when using only the two first strain invariants as features. Apart from accelerating learning and resulting in lossless dimensionality reduction, using invariants also results in a surrogate which is frame invariant under small strains. For comparison, we also train a conventional FNN on the same set of features, but those are unsurprisingly not enough to describe general strain states and much of the material response is interpreted by the FNN as observation noise. We zoom into the first part of the learning curves in \cref{FIelasticvsfnn2}, this time also showing single standard deviation uncertainty bands coming from the variance among the \SI{50}{} training datasets. The hybrid network outperforms conventional FNNs in the low data regime and tends to be less sensitive to changes in dataset starting from about \SI{20}{} training paths. Nevertheless, the extra flexibility of conventional FNNs allow them to achieve lower validation errors if significantly more training paths are used.

\begin{figure}
  \centering
  \begin{subfigure}[c]{0.33\textwidth}
    \includegraphics[scale=0.65]{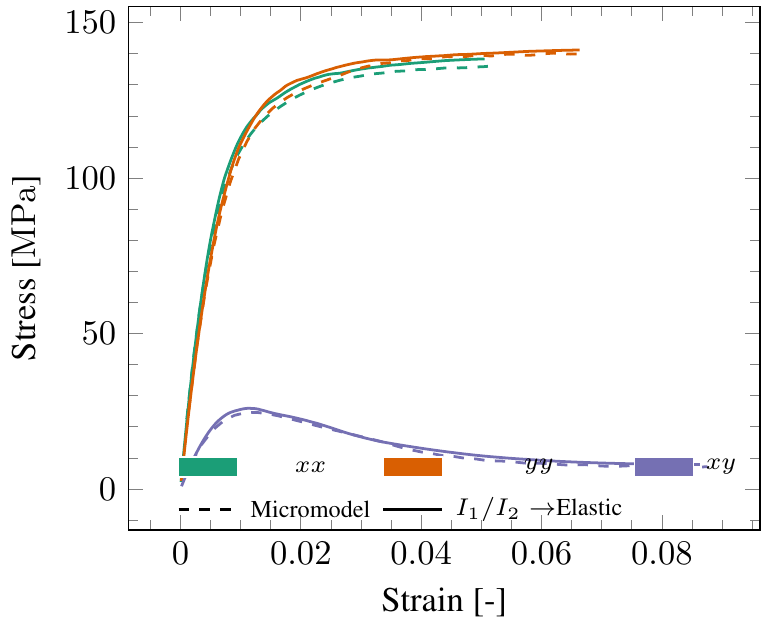}
    \caption{Monotonic}
  \end{subfigure}
  \begin{subfigure}[c]{0.33\textwidth}
    \includegraphics[scale=0.65]{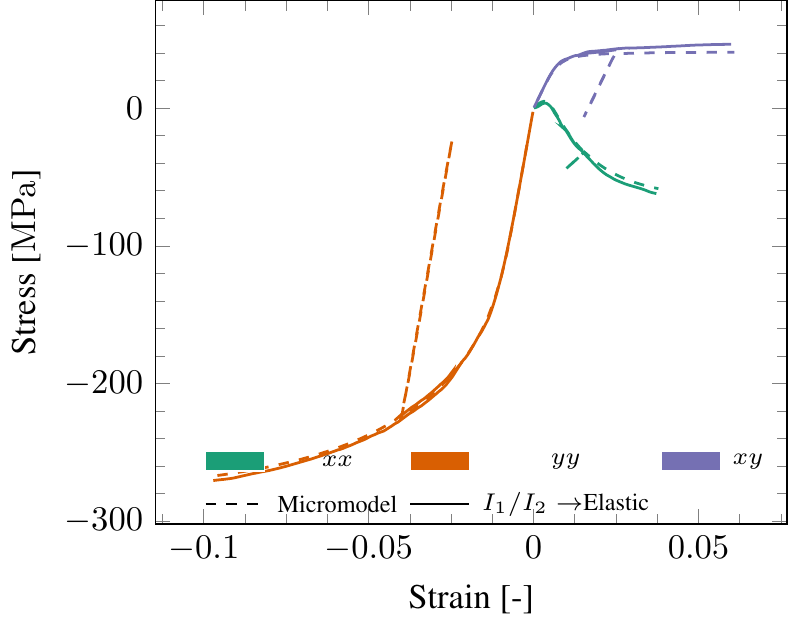}
    \caption{Unloading-reloading}
  \end{subfigure}
  \begin{subfigure}[c]{0.33\textwidth}
    \includegraphics[scale=0.65]{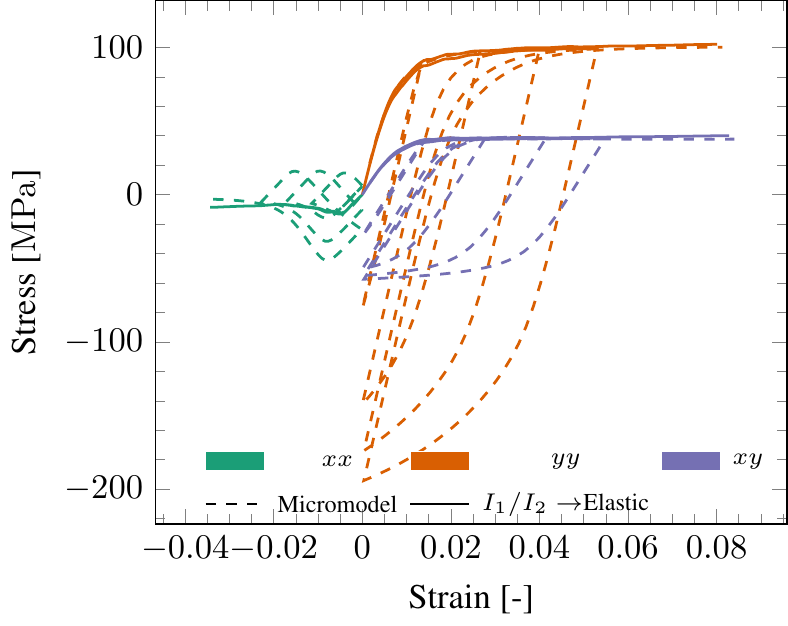}
    \caption{Slow cycling}
  \end{subfigure}
  \caption{Performance of the elastic decoder model for different test scenarios.}
  \label{FIelastici1i2tests}
\end{figure}

Training the invariant-based hybrid network with the complete dataset of \SI{1500}{} curves leads to surrogates with validation errors of about \SI{4}{\mega\pascal}, accurately representing the monotonic behavior of the original micromodel. \cref{FIelastici1i2tests} shows representative predictions of this model for paths from the test set. As expected, this surrogate with no internal variables is not capable of predicting non-monotonic strain paths, and effectively behaves like a hyperelastic material model just as the conventional FNN would. 

Nevertheless, the flexible and interpretable encoder-decoder architecture of \cref{FIapproach} allows for new creative approaches in feature selection. As a demonstration, we keep the trained network of \cref{FIelastici1i2tests} intact and only modify its feature extractor to introduce a simple path-dependent mechanism:
\begin{equation}
  \boldsymbol{\varphi}_T \equiv \begin{bmatrix} \overline{I}^\varepsilon_1 & \overline{I}^\varepsilon_2 \end{bmatrix}_T = \underset{0<t<T}{\mrm{argmax}}\left(\left(I_1^\varepsilon\right)^2_t + \left(J_2^\varepsilon\right)_t\right)
\end{equation}
\noindent which freezes the evolution of $\props$ if the path becomes non-monotonic. Note that the network does not need to be retrained and this modification can be employed exclusively when making online predictions, as the new features reduce to the original ones for the monotonic paths used for training.

\begin{figure}
  \centering
  \begin{subfigure}[c]{0.4\textwidth}
    \includegraphics[scale=0.8]{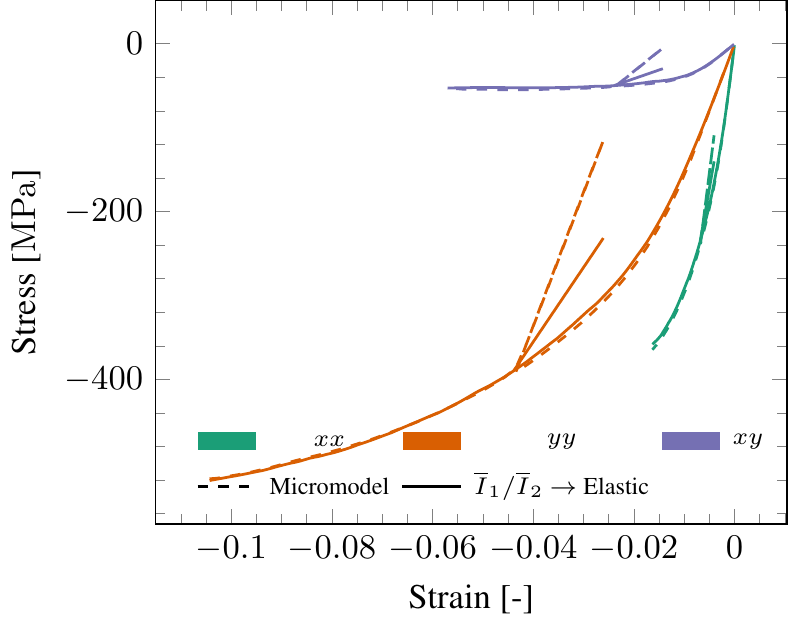}
  \end{subfigure}
  \begin{subfigure}[c]{0.4\textwidth}
    \includegraphics[scale=0.8]{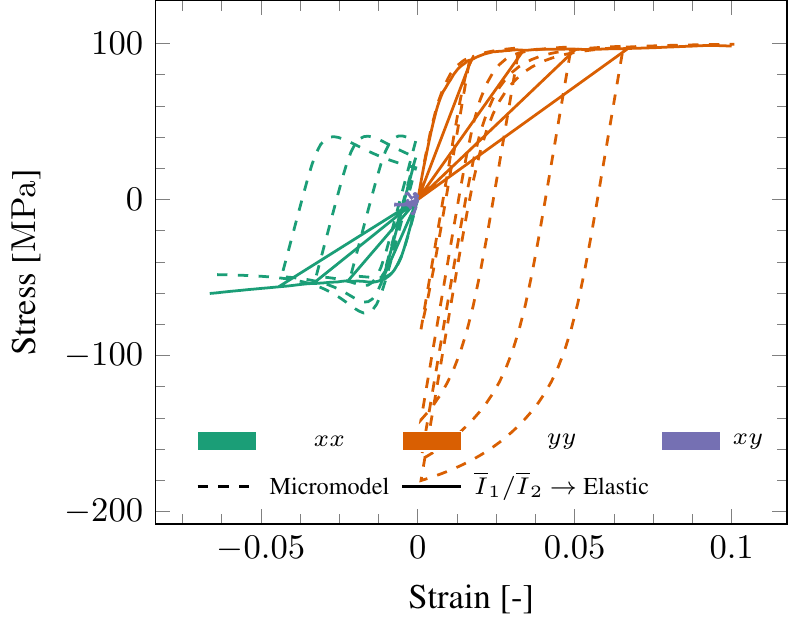}
  \end{subfigure}
  \caption{Predicting unloading with a linear-elastic decoder through history-aware feature extraction.}
  \label{FIelastici1i2mem}
\end{figure}

We plot two representative non-monotonic paths predicted by the modified model in \cref{FIelastici1i2mem}. From the hyperelastic behavior of \cref{FIelastici1i2tests}, the modified surrogate now behaves as a damage model: the non-linear material behavior is explained by a loss of stiffness which is made persistent by the history-aware feature extractor. Nevertheless, although an improvement to the original model, it is unreasonable to expect the physical bias introduced by a purely elastic model to reliably represent an elastoplastic micromodel. We therefore move to decoders with more relevant physics.

\subsection{$J_2$ decoder}

In this section we choose as decoder $\model$ the elastoplastic model of \cref{EQj2} with $J_2$ plastic flow. Standing on its own, the model is \textit{a priori} perfectly plastic (constant $\sigma_\mathrm{y}$)\footnote{We also experimented with models with linear hardening but since the encoder is in principle a universal approximator, no additional flexibility is obtained by assuming more complex hardening behavior. There are also no gains to be booked in terms of numerical stability, as a perfectly-plastic $J_2$ model is already unconditionally stable.}, but here we let its yield stress be controlled by the data-driven encoder: 
\begin{equation}
  \props = \begin{bmatrix} \sigma_\mrm{y} \end{bmatrix}
\end{equation}
\noindent while enforcing $10^1<\sigma_\mrm{y}<10^3$ and keeping the Young's modulus and Poisson's ratio fixed to values obtained from a single linear micromodel simulation. In contrast to the model with elastic decoder of the previous section, we now employ prior knowledge of the micromodel behavior and assume that all non-linearity should be explained by plasticity and do not let the elastic properties be dictated by the encoder. Still, the assumption of isotropic and incompressible plastic flow is a departure from the more complex pressure-dependent and non-associative behavior shown by the micromodel. Here we are therefore concerned with the effect of trading the flexibility of an elastic decoder for significantly more physical bias from a lower-fidelity representation of material behavior.

\begin{figure}
  \centering
  \includegraphics[scale=.9]{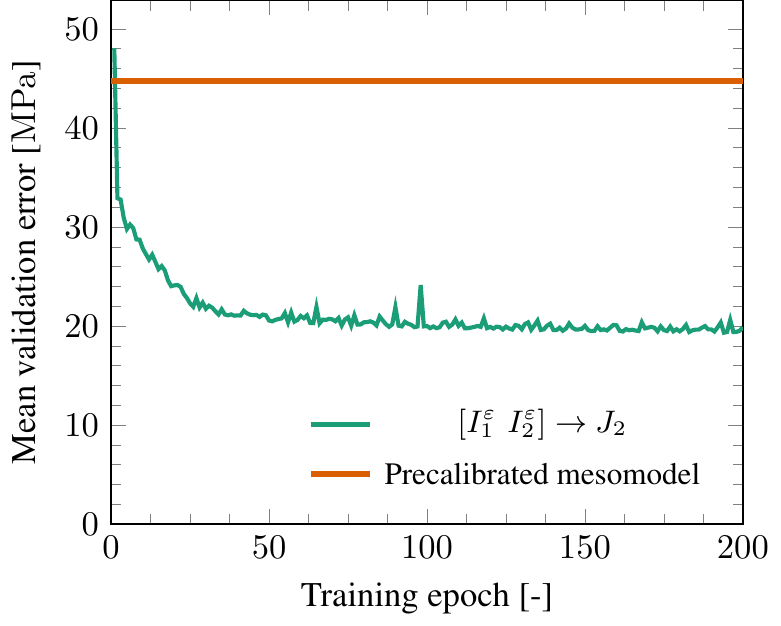}
  \caption{Evolution of the mean validation loss for the first 200 training epochs of a network with $J_2$ decoder. Single dataset with 1500 monotonic paths.}
  \label{FIj2loss}
\end{figure}

At this point it is interesting to compare the performance of the hybrid surrogate with predictions coming from the state-of-the-art mesoscale material model for polymer composites proposed by Vogler \etal \cite{voglerModelingInelasticDeformation2013}. It is an orthotropic elastoplastic model with pressure-dependent non-associative flow precalibrated with a small number of monotonic uniaxial and biaxial stress-strain curves obtained from simulations on the exact same micromodel of \cref{FImesh} (see \cite{vandermeerMicromechanicalValidationMesomodel2016} for details on the calibration procedure). For this section, we switch to a dataset in plane stress, allowing the $J_2$ model to describe richer nonlinear behavior under biaxial strain states.

\begin{figure}
  \centering
  \begin{subfigure}[c]{0.4\textwidth}
    \includegraphics[scale=0.8]{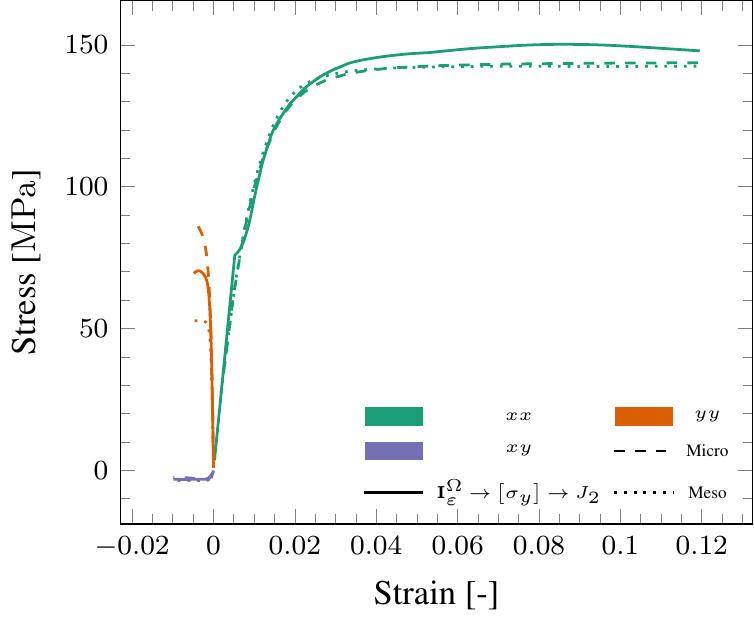}
    \caption{Predicting nonlinear hardening}
    \label{FIinvsigymonotonic1}
  \end{subfigure}
  \begin{subfigure}[c]{0.4\textwidth}
    \includegraphics[scale=0.8]{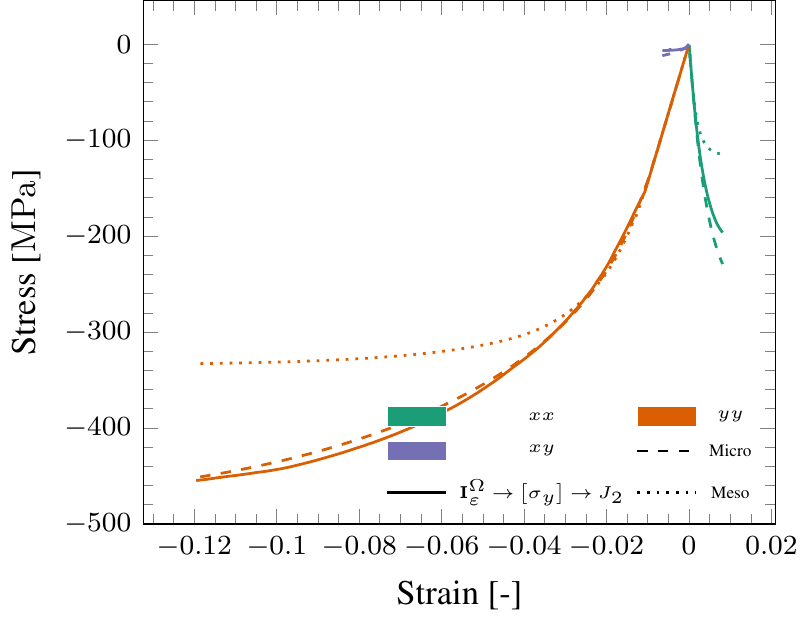}
    \caption{Predicting pressure dependency}
    \label{FIinvsigymonotonic2}
  \end{subfigure}
  \caption{Predictions from the network with $J_2$ decoder. Letting the yield stress evolve extends the model to more complex plasticity behavior.}
  \label{FIinvsigymonotonic}
\end{figure}

\cref{FIj2loss} shows the evolution of the validation set loss when training the $J_2$-decoded model with \SI{1500}{} plane stress training paths. The error quickly stabilizes at around \SI{20}{\mega\pascal}, significantly lower than the \SI{44}{\mega\pascal} average prediction error obtained with the precalibrated mesomodel. The added flexibility with respect to the original perfectly-plastic $J_2$ model can be seen in the test set curves plotted in \cref{FIinvsigymonotonic}: the data-driven encoder leads to correct predictions of nonlinear hardening (\cref{FIinvsigymonotonic1}) and pressure-dependent plastic flow (\cref{FIinvsigymonotonic2}). The figures also highlight the inability of the mesomodel to predict the behavior in certain regions of the strain space, particularly under compression-dominated scenarios.

The minimum validation error attained by the model is, however, nevertheless significantly higher than the \SI{4}{\mega\pascal} obtained with the elastic decoder of the previous section. This result is not entirely surprising, as the elastic decoder introduces much less bias into the model and therefore allows for a greater degree of flexibility when fitting monotonic data. On the other hand, what cannot be directly gleaned from \cref{FIj2loss} is that the $J_2$ decoder benefits from having physics-based memory coming from its internal variables that allows for making predictions of non-monotonic behavior based solely on our assumption that nonlinearities come from plastic strain and therefore without ever having to see it during training.

\begin{figure}
  \centering
  \begin{subfigure}[c]{0.4\textwidth}
    \includegraphics[scale=0.8]{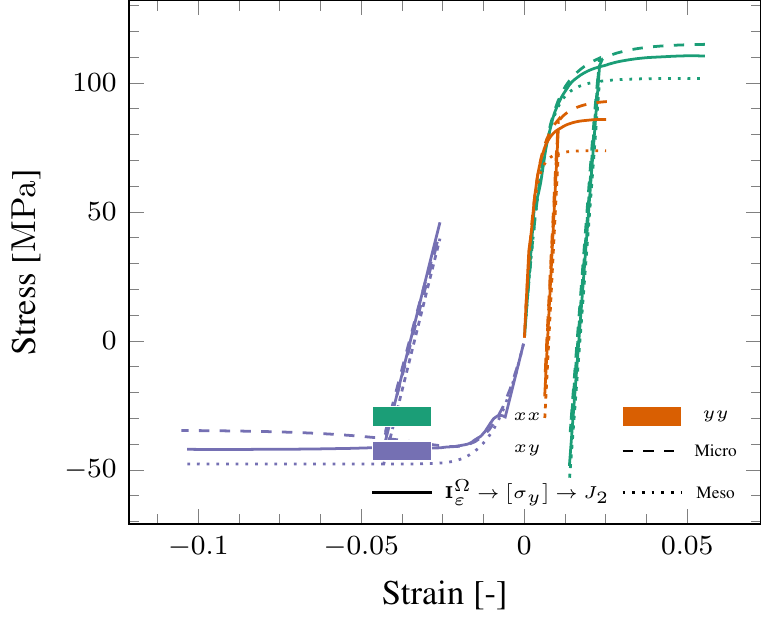}
  \end{subfigure}
  \begin{subfigure}[c]{0.4\textwidth}
    \includegraphics[scale=0.8]{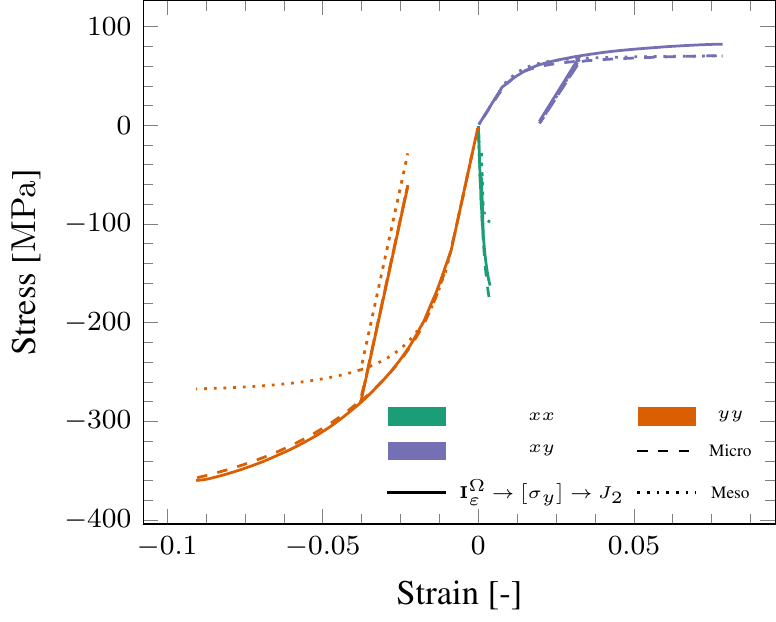}
  \end{subfigure}
  \caption{Network predictions with $J_2$ decoder for unloading paths after being trained exclusively with monotonic paths.}
  \label{FIinvsigyunloading}
\end{figure}

\begin{figure}
  \centering
  \begin{subfigure}[c]{0.33\textwidth}
    \includegraphics[scale=0.65]{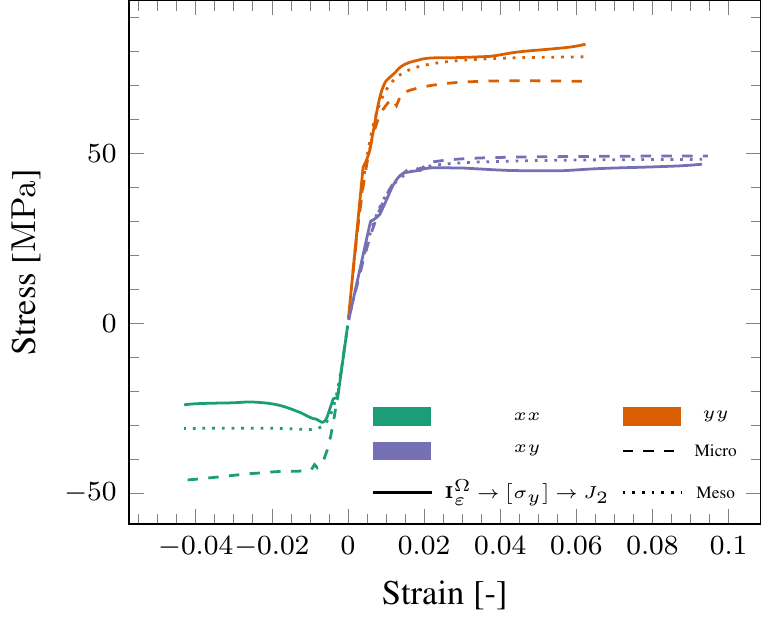}
  \end{subfigure}
  \begin{subfigure}[c]{0.33\textwidth}
    \includegraphics[scale=0.65]{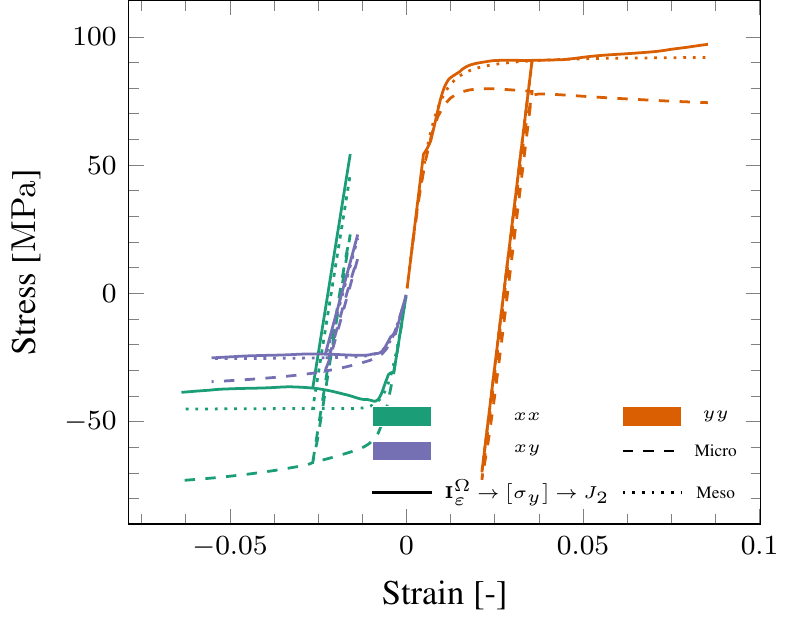}
  \end{subfigure}
  \begin{subfigure}[c]{0.33\textwidth}
    \includegraphics[scale=0.65]{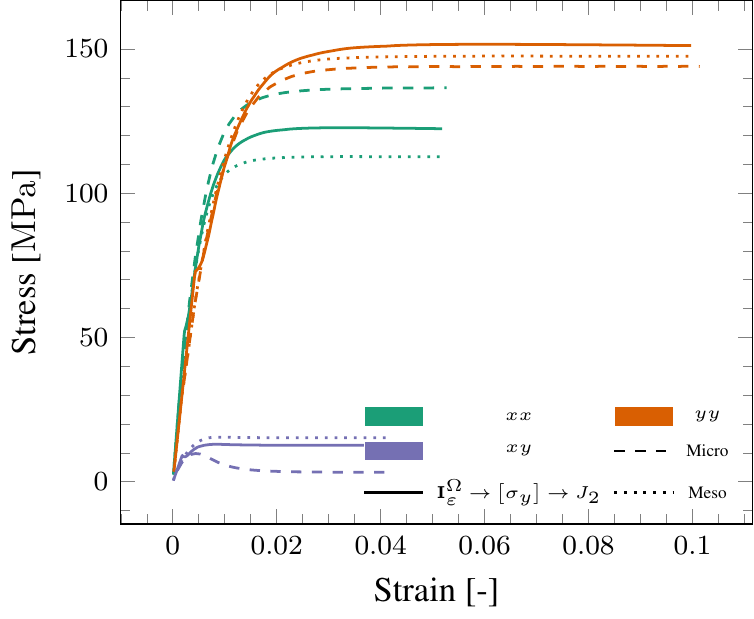}
  \end{subfigure}
  \caption{Examples of strain paths not well predicted by the $J_2$ decoded model.}
  \label{FIinvsigymonobadpredictions}
\end{figure}

In \cref{FIinvsigyunloading} we plot predictions of the $J_2$ surrogate for two different unloading-reloading paths from the test dataset. The model predicts unloading very well without being trained for it. Nevertheless, as \cref{FIj2loss} suggests, the model struggles to predict monotonic behavior under a number of different scenarios, from which it follows that any non-monotonic predictions along the same directions will also be inaccurate. \cref{FIinvsigymonobadpredictions} shows three examples of this behavior. 

The choice of decoder therefore involves a tradeoff between bias and flexibility that can be deceiving to base solely on validation error computed on monotonic data. Indeed, the decoder used in the next section outperforms $J_2$-based decoders in most situations, but nevertheless a choice for the simpler decoder might still be justified --- \eg if the unconditional numerical stability of a $J_2$ decoder is desirable.

\subsection{Non-associative pressure-dependent elastoplastic decoder}

As one final exploration on model selection, we use as decoder the same elastoplastic model by Melro \etal\ used to describe the matrix material at the microscale \cite{melro_micromechanical_2013}. As mentioned in \cref{SEdecoders}, this model is the natural choice for \model, as it attempts to explain the observed microscopic non-linear behavior with the same model from which the behavior arises. As before we keep the elastic properties of the model intact and let only the yield stresses and the plastic Poisson's ratio change in time:
\begin{equation}
  \props = \begin{bmatrix}\sigma_\tens & \frac{\sigma_\comp}{\sigma_\tens} & \nu_\p \end{bmatrix}
\end{equation}
\noindent where $10^1<\sigma_\tens<10^4$, $0<\nu_\p<0.5$ and $1<\frac{\sigma_\comp}{\sigma_\tens}<100$. We opt for the ratio $\frac{\sigma_\comp}{\sigma_\tens}$ instead of simply $\sigma_\comp$ in order to also enforce $\sigma_\comp>\sigma_\tens$.

\begin{figure}
  \centering
  \includegraphics[scale=1]{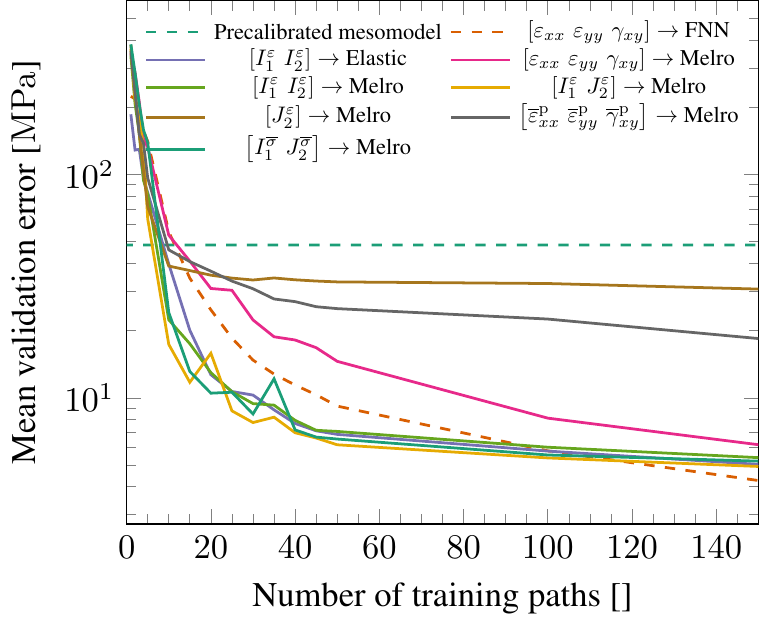}
  \caption{Expected validation errors for Melro-decoded surrogates with different feature extractors (averages over \SI{50}{} datasets).}
  \label{FImelrovali}
\end{figure}

We expand upon the feature selection study of \cref{FIelasticvsfnn} by looking at several feature extractors coming both directly from strains and from the output of a precalibrated Melro model $\overline\model$ with the same properties used at the microscale (\cref{FIarch2}). Aside from the familiar choice of strain features (\textit{$\left[\varepsilon_{xx}\,\,\varepsilon_{yy}\,\,\gamma_{xy}\right]\rightarrow$ Melro}), we look into invariants of the strain tensor (\textit{$\left[I_1^\varepsilon\,\,I_2^\varepsilon\right]\rightarrow$ Melro}), combinations including invariants of the deviatoric strain tensor (\textit{$\left[J_2^\varepsilon\right]\rightarrow$ Melro}, \textit{$\left[I_1^\varepsilon\,\,J_2^\varepsilon\right]\rightarrow$ Melro}), plastic strain internal variables coming from the precalibrated feature extractor (\textit{$\left[\overline{\varepsilon}^\mathrm{p}_{xx}\,\,\overline{\varepsilon}^\mathrm{p}_{yy}\,\,\overline{\gamma}^\mathrm{p}_{xy}\right]\rightarrow$ Melro}) and stress invariants coming from the extractor (\textit{$\left[I_1^{\overline{\sigma}}\,\,J_2^{\overline{\sigma}}\right]\rightarrow$ Melro}). We also include the precalibrated mesomodel by Vogler \etal \cite{voglerModelingInelasticDeformation2013} and selected curves from \cref{FIelasticvsfnn1} for comparison purposes. As before, we train \SI{50}{} networks of each type for each size of dataset ranging from \SI{1}{} to \SI{150}{} paths drawn from the original dataset with \SI{1500}{} paths. Each trained network is then used to compute the validation error over the \SI{500}{} monotonic validation paths and the \SI{150}{} test paths (\SI{50}{} extra monotonic paths, \SI{50}{} paths with unloading-reloading and \SI{50}{} slow cycle paths). This results in an extensive study comprising \SI{6800}{} trained networks and over one million test set simulations.

\begin{figure}
  \centering
  \begin{subfigure}[c]{0.4\textwidth}
    \includegraphics[scale=0.8]{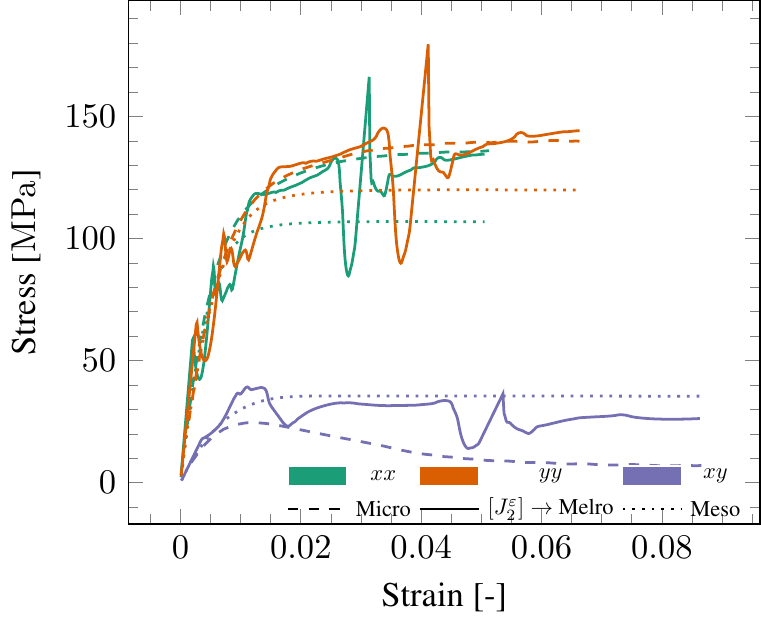}
    \caption{Model with second deviatoric strain invariant as feature}
  \end{subfigure}
  \begin{subfigure}[c]{0.4\textwidth}
    \includegraphics[scale=0.8]{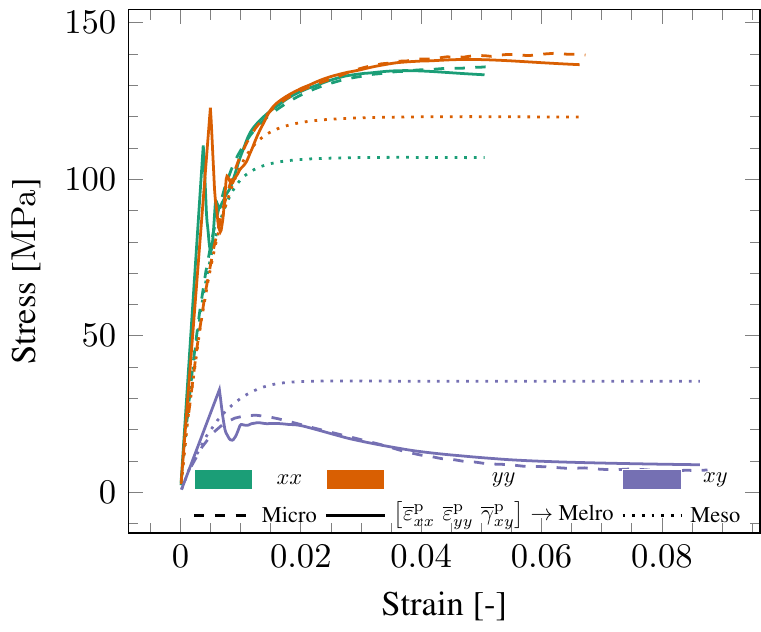}
    \caption{Model with plastic strain features}
  \end{subfigure}
  \caption{Monotonic test set predictions from feature-deficient Melro models (complete training dataset with 1500 paths).}
  \label{FImelroepspj2}
\end{figure}

Results are summarized in \cref{FImelrovali}, with each point in a curve being the average over \SI{50}{} networks. Once again using invariants as features proves beneficial, leading to lossless dimensionality reduction and frame invariant surrogates. All tested models perform better than the precalibrated mesomodel, with a gap of more than one order of magnitude for the best performing surrogates. Interestingly, models with Melro-based decoders seem to learn as fast and be as flexible as models with elastic decoders, already for the monotonic curves in the validation dataset. This suggests that the new decoder does not impose extra undesirable bias in learning the specific material behavior treated here other than the assumptions that had already been introduced by elasticity (\eg symmetries and couplings encoded by the elastic stiffness tensor). Any benefits reaped when extrapolating to non-monotonic paths, as we will see in the following, are therefore obtained at a negligible price in terms of monotonic behavior accuracy. This stands in contrast with the discussion on the $J_2$ decoder of the previous section.

\begin{figure}
  \centering
  \begin{subfigure}[c]{0.4\textwidth}
    \includegraphics[scale=0.8]{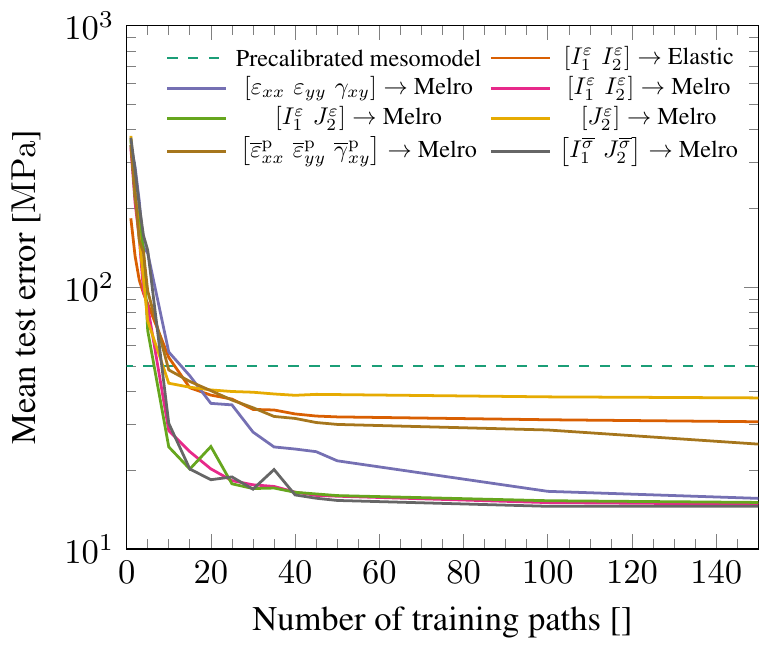}
    \caption{Errors for complete paths}
    \label{FImelrounloading1}
  \end{subfigure}
  \begin{subfigure}[c]{0.4\textwidth}
    \includegraphics[scale=0.8]{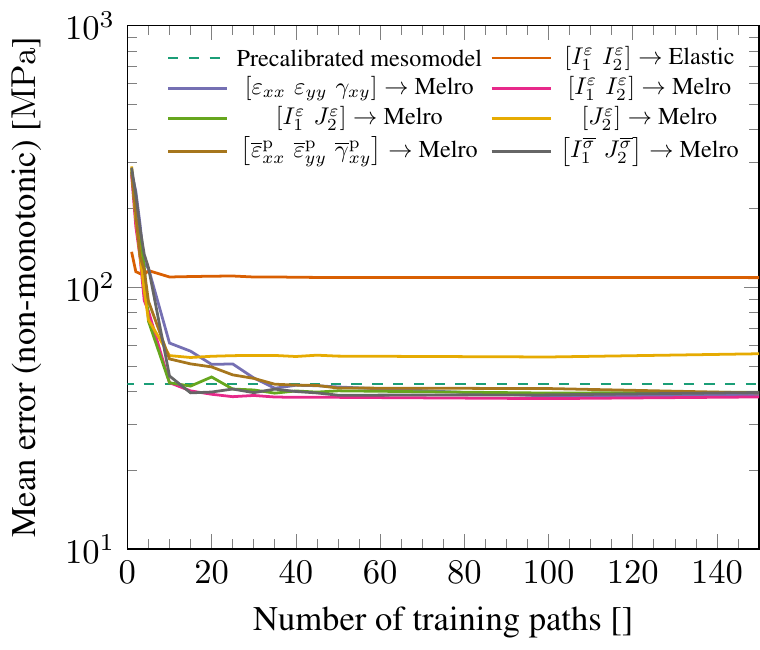}
    \caption{Errors including only unloading/reloading steps}
    \label{FImelrounloading2}
  \end{subfigure}
  \caption{Learning curves for unloading-reloading test errors of Melro-decoded surrogates (averages of \SI{50}{} datasets).}
  \label{FImelrounloading}
\end{figure}

Although \cref{FImelrovali} is not enough to discern between several of our encoder choices, it is interesting to take a closer look at the two clearly underperforming options. \cref{FImelroepspj2} shows predictions from \textit{$\left[J_2^\varepsilon\right]\rightarrow$ Melro} and \textit{$\left[\overline{\varepsilon}^\mathrm{p}_{xx}\,\,\overline{\varepsilon}^\mathrm{p}_{yy}\,\,\overline{\gamma}^\mathrm{p}_{xy}\right]\rightarrow$ Melro} for the same monotonic test path. The model with a single feature struggles to predict the entirety of the path, indicating that further reducing the dimensionality of the feature space is not possible for this dataset. The oscillatory stress predictions make this model unsuitable for online stress evaluation in a multiscale setting. For the model with plastic strain features, the feature extractor shows no plastic strains until high stress levels while in the micromodel plasticity starts much earlier, forcing the surrogate to remain in the elastic regime until a sudden jump brings it back to the expected path.

Moving to unloading-reloading paths, we compare the performance of different feature sets by plotting the average test error over the \SI{50}{} unloading-reloading paths in \cref{FImelrounloading1}. Here an interesting observation can be made: even the surrogate \textit{$\left[I^\varepsilon_1\,\,I^\varepsilon_2\right]\rightarrow$ Elastic} --- which cannot predict unloading at all --- attains a lower test error than the precalibrated mesomodel. This apparent contradiction can be explained by plotting in \cref{FImelrounloading2} the average error computed only at unloading or reloading time steps: use of an elastic decoder --- and therefore of a conventional FNN or an RNN trained with insufficient data --- excels at predicting monotonic response but is consistently inaccurate for non-monotonic paths and shows little improvement when more monotonic paths are added to the training dataset. 

\begin{figure}
  \centering
  \begin{subfigure}[c]{0.4\textwidth}
    \includegraphics[scale=0.80]{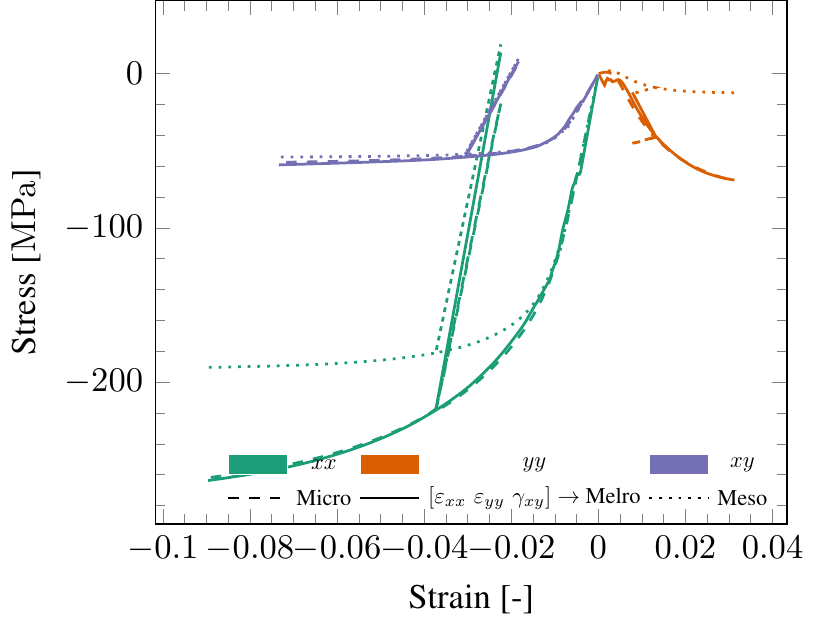}
  \end{subfigure}
  \begin{subfigure}[c]{0.4\textwidth}
    \includegraphics[scale=0.80]{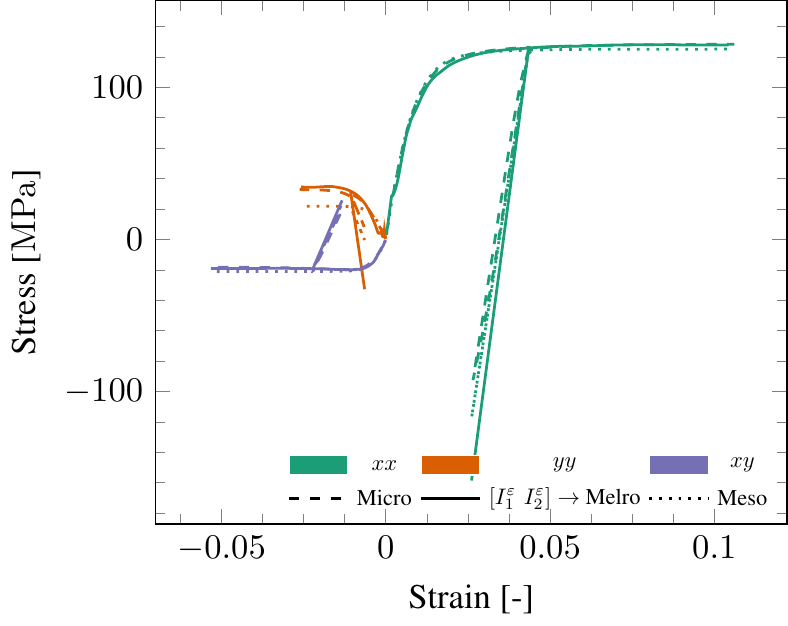}
  \end{subfigure}\\
  \begin{subfigure}[c]{0.4\textwidth}
    \includegraphics[scale=0.80]{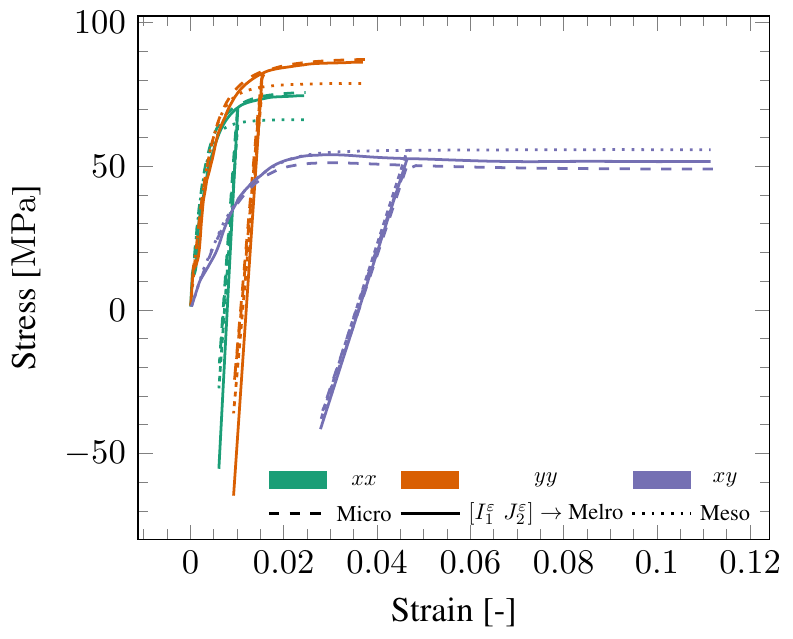}
  \end{subfigure}
  \begin{subfigure}[c]{0.4\textwidth}
    \includegraphics[scale=0.80]{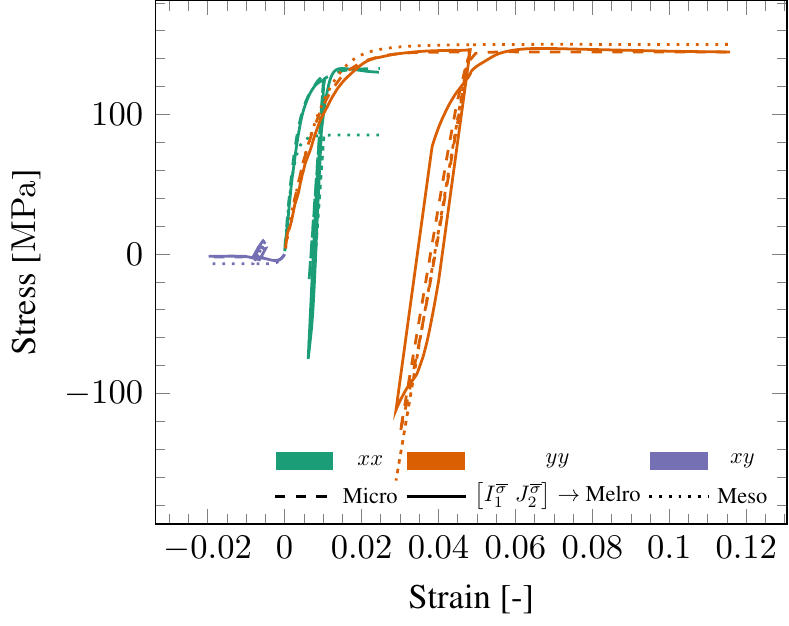}
  \end{subfigure}
  \caption{Response of Melro-decoded surrogates with different features for selected unloading/reloading test paths (\SI{1500}{} monotonic training paths).}
  \label{FImelrounloadingexamples}
\end{figure}

In contrast, the best-performing Melro models are consistently more accurate than the precalibrated mesomodel even when trained on very little data. We plot in \cref{FImelrounloadingexamples} selected representative unloading paths from the test dataset for four of the surrogates. Unloading is once again well captured without having been seen during training, and since it emerges from a purely physical mechanism, it is reasonable to expect unloading at different points along the path to yield comparable results (\cf \cref{FIlstmunloading}). 
Nevertheless, relatively small differences in unloading slope can still lead to large differences in stress at the end of the unloading branches. Furthermore, the model can struggle with tension-compression switches and predict spurious hysteresis loops.

\begin{figure}
  \centering
  \begin{subfigure}[c]{0.4\textwidth}
    \includegraphics[scale=0.8]{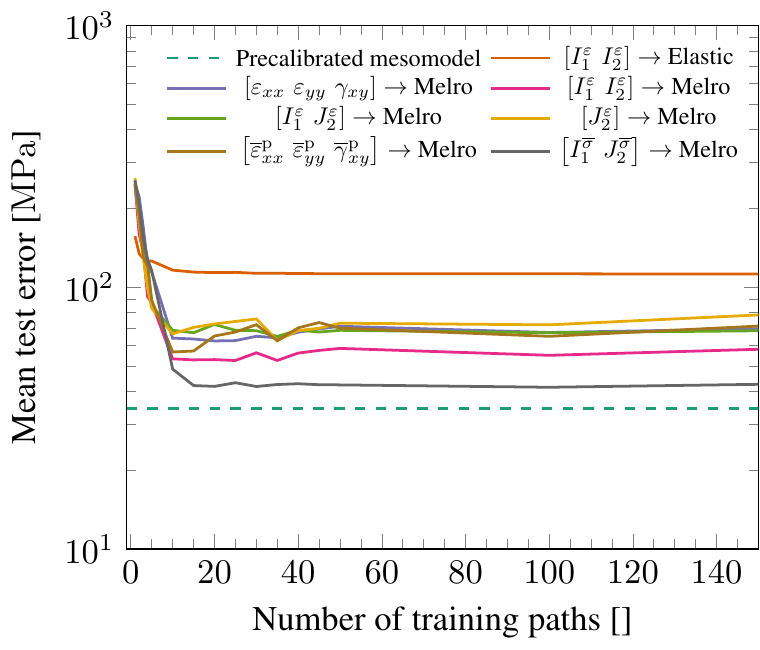}
    \caption{Errors for complete paths}
  \end{subfigure}
  \begin{subfigure}[c]{0.4\textwidth}
    \includegraphics[scale=0.8]{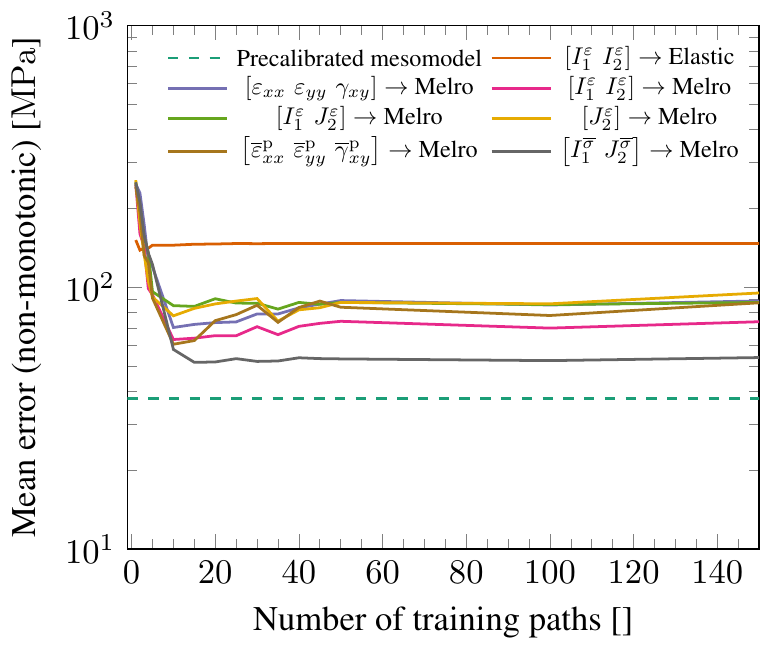}
    \caption{Errors including only unloading/reloading steps}
  \end{subfigure}
  \caption{Slow cycling test errors for Melro-decoded surrogates (averages of \SI{50}{} datasets for each size).}
  \label{FImelrocycles}
\end{figure}

Indeed, we observe a consistent inability by the models to properly predict switches between tension and compression within the same path. This becomes clear when looking at slow cycling test paths composed of several of these switches (\cref{FItestexamplescycles}). We plot learning curves for the test error on slow cycling paths in \cref{FImelrocycles}, for complete paths as well as exclusively for the non-monotonic branches of the paths. In contrast with results up until now, here we see larger differences in performance for different feature sets. As expected, elastic decoders are once again shown to be unsuitable to predict non-monotonic paths, and the difference here is even more pronounced than in for single-unloading paths (\cf \cref{FImelrounloading}) as most of the path is composed of unloading/reloading branches. The model encoded with stress invariants coming from an elastoplastic feature extractor performs best among the models we test. But crucially, none of the surrogates manages to surpass the precalibrated mesomodel in this case. 

\begin{figure}
  \centering
  \begin{subfigure}[c]{0.4\textwidth}
    \includegraphics[scale=0.8]{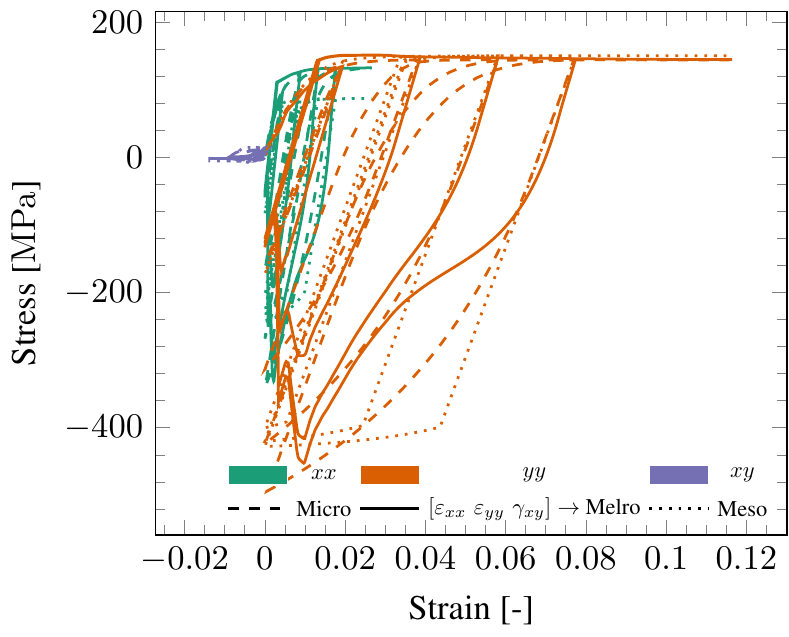}
  \end{subfigure}
  \begin{subfigure}[c]{0.4\textwidth}
    \includegraphics[scale=0.8]{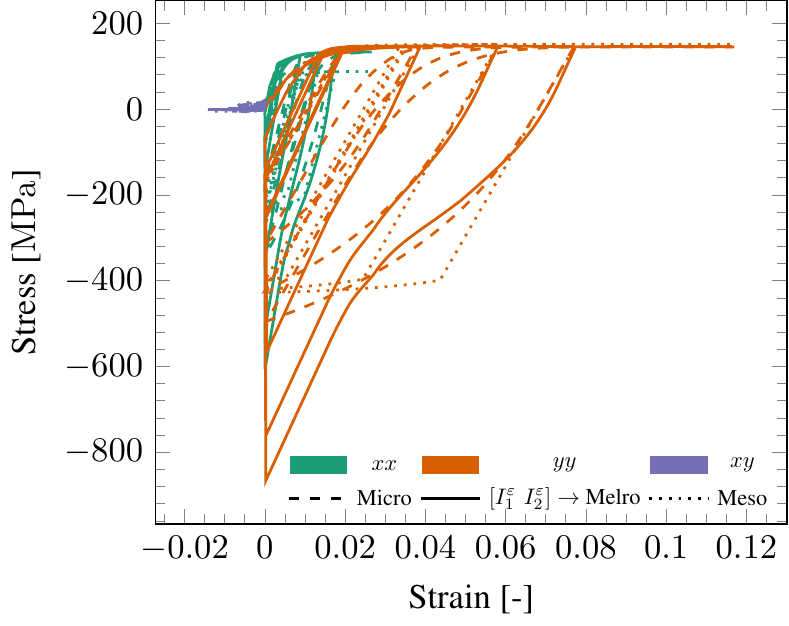}
  \end{subfigure}\\
  \begin{subfigure}[c]{0.4\textwidth}
    \includegraphics[scale=0.8]{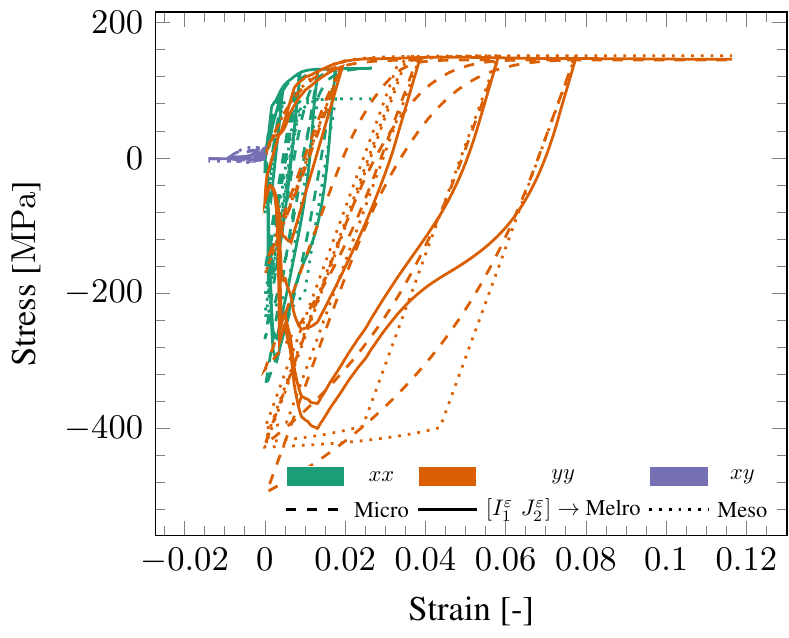}
  \end{subfigure}
  \begin{subfigure}[c]{0.4\textwidth}
    \includegraphics[scale=0.8]{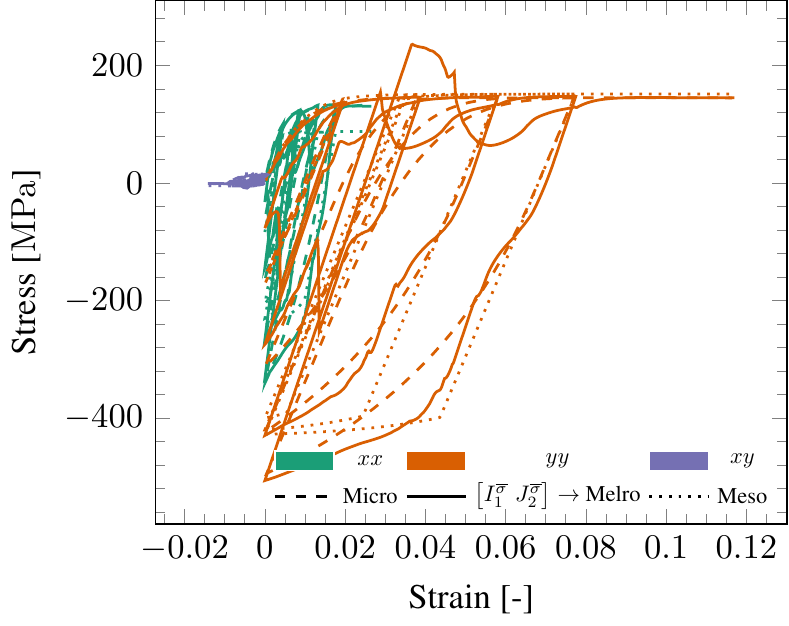}
  \end{subfigure}
  \caption{Response of Melro-decoded surrogates with different features for selected slow cycling test paths (\SI{1500}{} monotonic training paths).}
  \label{FImelrocyclescomp}
\end{figure}

As a demonstration, we select a representative path from the test dataset and plot predictions made with four different feature sets in \cref{FImelrocyclescomp}. 
As expected, larger errors are observed for more pronounced tension-compression switches as models either over- or undershoot the stress levels at compression-tension switch points. Interestingly, most models manage to converge back to the correct stress path after reloading, since hardening behavior is completely dictated by their non-recurrent data-driven encoders. The exception is the model with stress invariant features (\textit{$\left[I_1^{\overline{\sigma}}\,\,J_2^{\overline{\sigma}}\right]\rightarrow$ Melro}), performing significantly better than the rest 
but showing a number of undesired oscillations in stress response due to the (physically) recurrent nature of its features forcing its neural network encoder to operate in extrapolation.

\subsection{\fetwo\ example}

\begin{figure}
  \centering
  \includegraphics[scale=2]{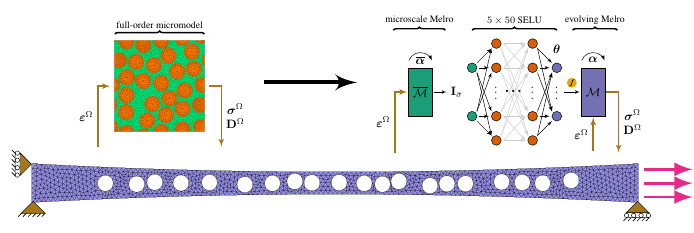}
  \caption{\fetwo\ example: geometry, mesh and boundary conditions. Full-order (left) and surrogate-based (right) \fetwo\ simulations are compared.}
  \label{FIfe2}
\end{figure}

We conclude our discussion with an \fetwo\ demonstration using the proposed hybrid surrogate. We model the tapered macroscopic bar with geometry and boundary conditions shown in \cref{FIfe2}. The model is meshed with 1620 linear triangles with a single Gauss point each and is loaded in tension until plastic strain localization takes place. The combination of the tapered geometry with the several circular voids along the model result in a complex range of stress states throughout the model. In contrast to the cases considered so far, this example also covers non-proportional strain paths. To facilitate convergence, the substepping approach proposed in \cite{somerSubsteppingSchemeMultiscale2009} is employed and an adaptive stepping algorithm is used at the macroscale that automatically reduces time step size and recomputes the current increment if either the micro- or macroscopic Newton-Raphson solver fails to converge.

\begin{figure}
  \centering
  \begin{subfigure}[c]{0.4\textwidth}
    \includegraphics[scale=.8]{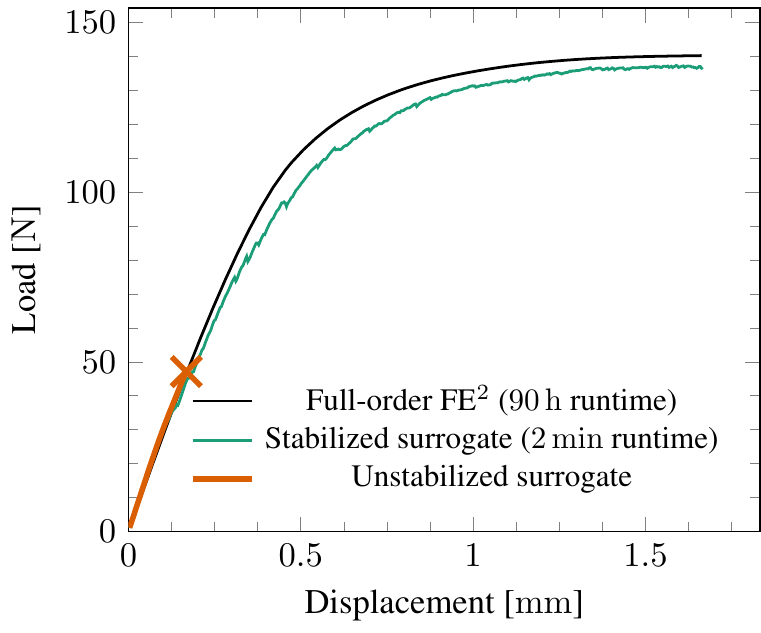}
    \caption{Monotonic}
    \label{FIholedogmnnvsfull1}
  \end{subfigure}
  \begin{subfigure}[c]{0.4\textwidth}
    \includegraphics[scale=.8]{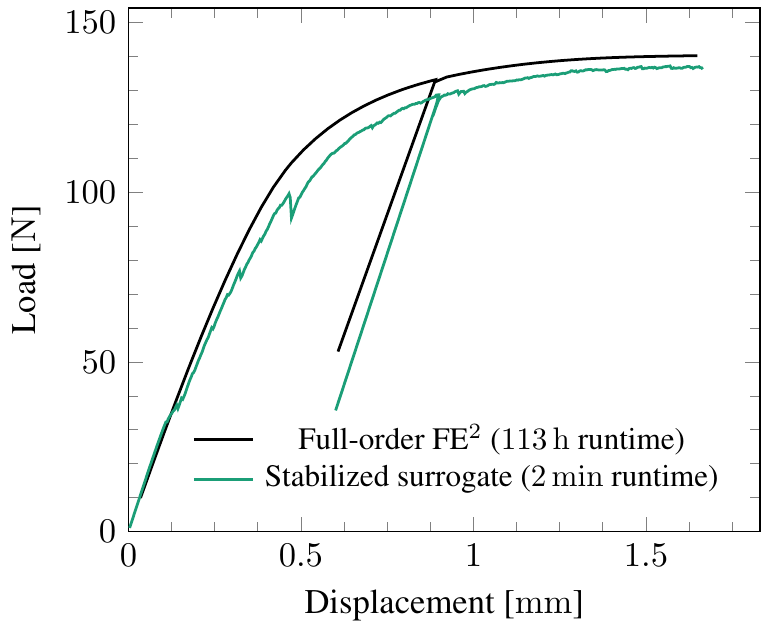}
    \caption{Unloading/reloading}
    \label{FIholedogmnnvsfull2}
  \end{subfigure}
  \caption{\fetwo example: load-displacement curves with and without online stabilization, compared to the ground-truth solution.}
  \label{FIholedogmnnvsfull}
\end{figure}

We use the \textit{$\left[I_1^{\overline{\sigma}}\,\,J_2^{\overline{\sigma}}\right]\rightarrow$ Melro} model of the previous section as surrogate, trained on the complete set of \SI{1500}{} monotonic training strain paths. The global load-displacement curve at the right edge of the model is plotted for the full-order \fetwo\ solution and using the hybrid surrogate in \cref{FIholedogmnnvsfull1}. Since we update decoder properties in an explicit fashion (\ie once per time step, see \cref{ALwrapper}), we use a displacement increment $\Delta u = \SI{3.5e-3}{\mm}$ for the approximate model, \SI{10}{} times smaller than the one used for the full-order model.

As mentioned in \cref{SEstabilization}, the model by Melro \etal can suffer from numerical stability issues even with fixed material properties, and it is reasonable to expect these issues to become worse when letting properties evolve with time. Indeed, with no additional stabilization the model using the network fails to converge at the point marked in \cref{FIholedogmnnvsfull1}. In contrast, the stabilization procedure of \cref{SEstabilization} allows for a complete path to be obtained. For this first result, we stabilize the network for 5 epochs with a learning rate of \SI{1e-5}{} for the stabilization loss (\cref{EQstabilityloss}) and \SI{1e-9}{} for retraining on a single monotonic training path selected at random. We also consider a model with an unloading/reloading switch after the onset of macroscopic plasticity. Results are shown in \cref{FIholedogmnnvsfull2}. The surrogate approximates the full-order behavior fairly accurately and several orders of magnitude faster than the full-order model. 

\begin{figure}
  \centering
  \begin{subfigure}[c]{0.4\textwidth}
    \includegraphics[scale=0.8]{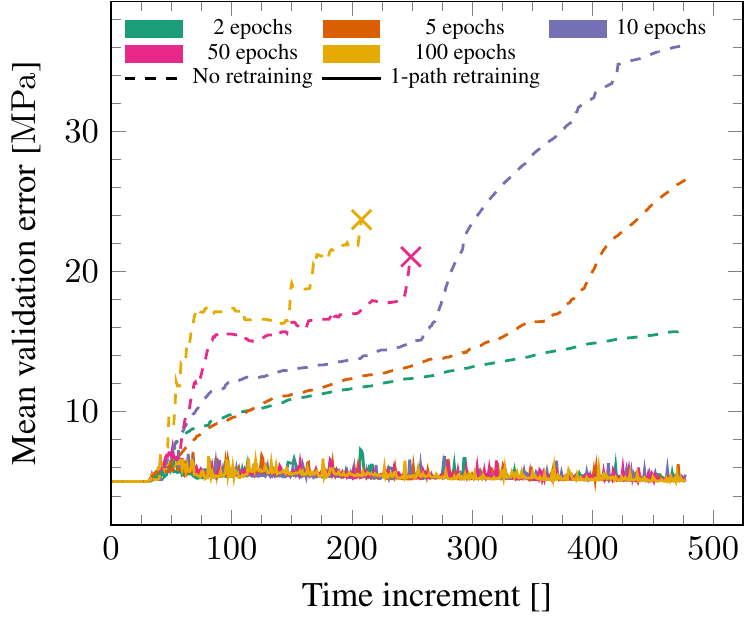}
    \caption{\SI{500}{}-path validation loss after stabilization}
    \label{FIholedogvalislodis1}
  \end{subfigure}
  \begin{subfigure}[c]{0.4\textwidth}
    \includegraphics[scale=0.8]{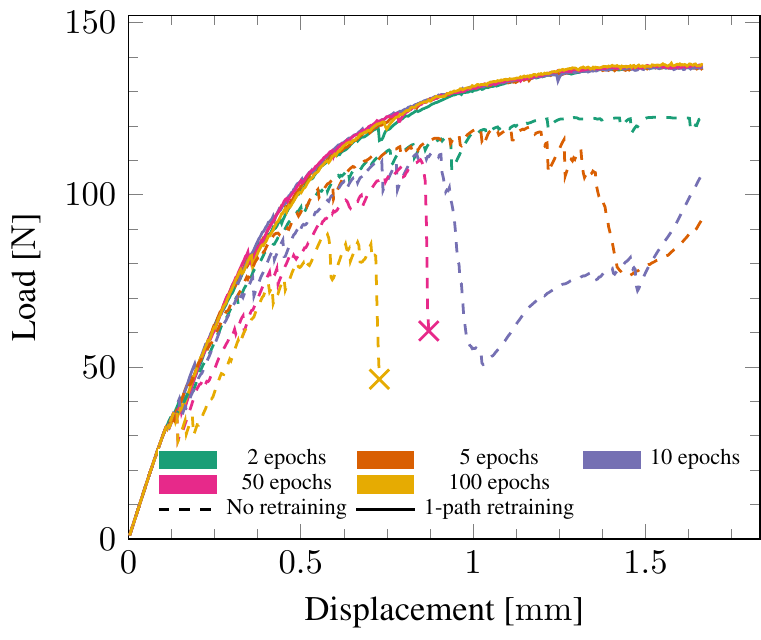}
    \caption{Load-displacement curves}
    \label{FIholedogvalislodis2}
  \end{subfigure}
  \caption{Performance of the surrogate model for stabilization strategies of varying intensities with and without retraining after stabilization.}
  \label{FIholedogvalislodis}
\end{figure}

We now look closer on the performance of the proposed online stabilization approach. We empirically find that retraining the network until every violating material point is fully stabilized is not strictly necessary in order to achieve convergence, and therefore opting for a small number of stabilization epochs proves to be an efficient approach. It is nevertheless interesting to investigate the impact of the number of stabilization epochs and of the subsequent retraining minibatch on the original dataset. We solve the monotonic example of \cref{FIholedogmnnvsfull1} with different numbers of stabilization/retraining epochs ranging from \SI{2}{} to \SI{100}{} and compute the validation loss (on the \SI{500}{}-path validation set used for model selection) at the end of every macroscopic time increment in order to keep track of how much the stabilized network deviates from its original pretrained state. 

Results are shown together with the corresponding load-displacement curves in \cref{FIholedogvalislodis}. All curves remain stable at first, as stabilization is only triggered when the first unstable points are detected. From that point, models which do not undergo retraining after stabilization lose accuracy at a rate proportional to the number of stabilization epochs. However, this unintuitively does not lead to improved global stability: the loss of accuracy by the surrogate leads to spurious global softening (\cf \cref{FIholedogvalislodis2}) which in turn leads to further need for stabilization. Models stabilized for \SI{50}{} and \SI{100}{} epochs continuously fail to converge and we opt for terminating the simulation after \SI{100}{} cancelled time increments. On the other hand, models retrained with as little as a single strain path (out of the original \SI{1500}{}) after each stabilization epoch are able to maintain the original model accuracy while offering enough stability gains to allow the simulation to converge until the final step, with little change in global behavior for different stabilization regimes.

\begin{figure}
  \centering
  \begin{subfigure}[c]{0.4\textwidth}
    \includegraphics[scale=0.8]{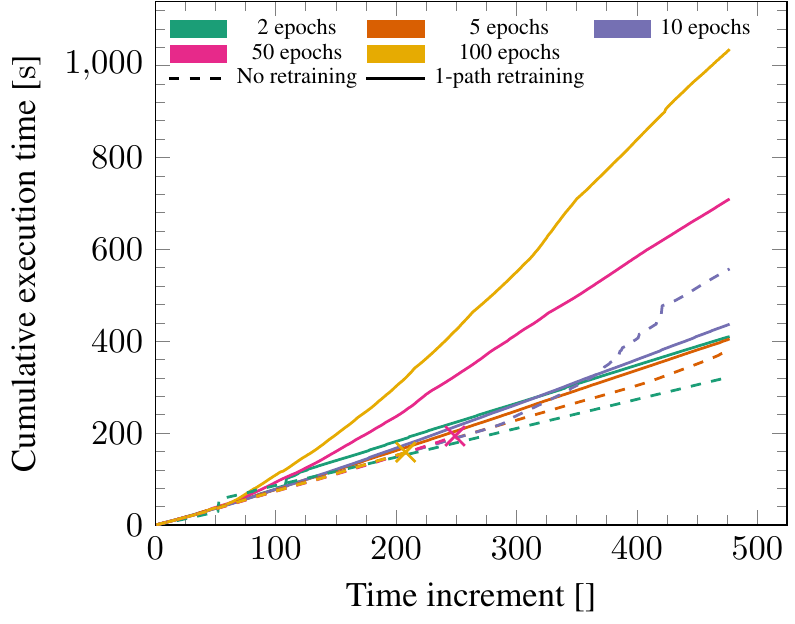}
    \caption{Computational overhead due to stabilization}
    \label{FIholedogexectimesstabs1}
  \end{subfigure}
  \begin{subfigure}[c]{0.4\textwidth}
    \includegraphics[scale=0.8]{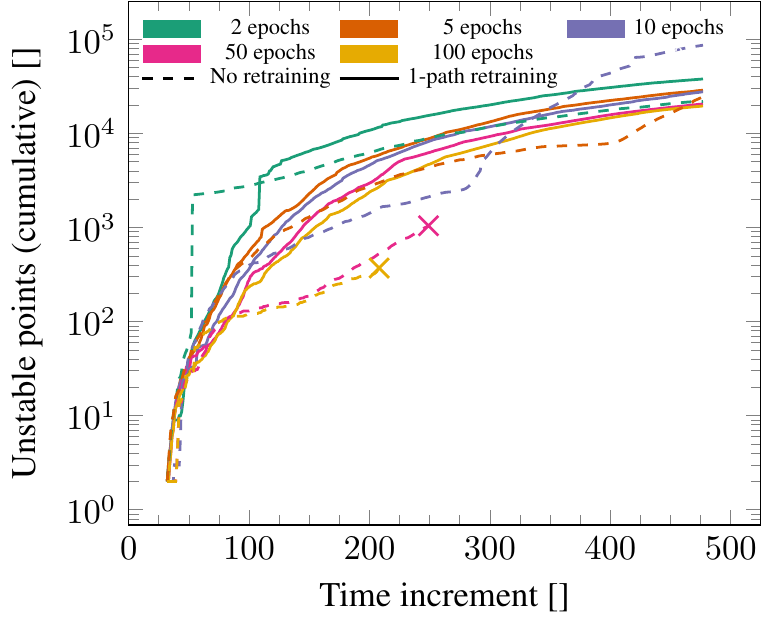}
    \caption{Number of detected unstable strain states}
    \label{FIholedogexectimesstabs2}
  \end{subfigure}
  \caption{Impact of stabilization regime on execution time and number of unstable points throughout the simulation.}
  \label{FIholedogexectimesstabs}
\end{figure}

More insight can be obtained on the different stabilization strategies by plotting the cumulative execution time of the simulation and the cumulative number of detected unstable strain states with time increments for different numbers of stabilization epochs. Results can be seen in \cref{FIholedogexectimesstabs}. In general, simulations without retraining tend to run faster and result in improved stability, although any gains are quickly overshadowed by losses in accuracy (\cf \cref{FIholedogvalislodis}). Stabilizing for more epochs results in a reduction in the total number of unstable points detected, but beyond 5 epochs this does not result in an overall reduction in the computational cost of the simulation given the increased effort spent on individual stabilization operations.

\begin{figure}
  \centering
  \includegraphics[scale=1]{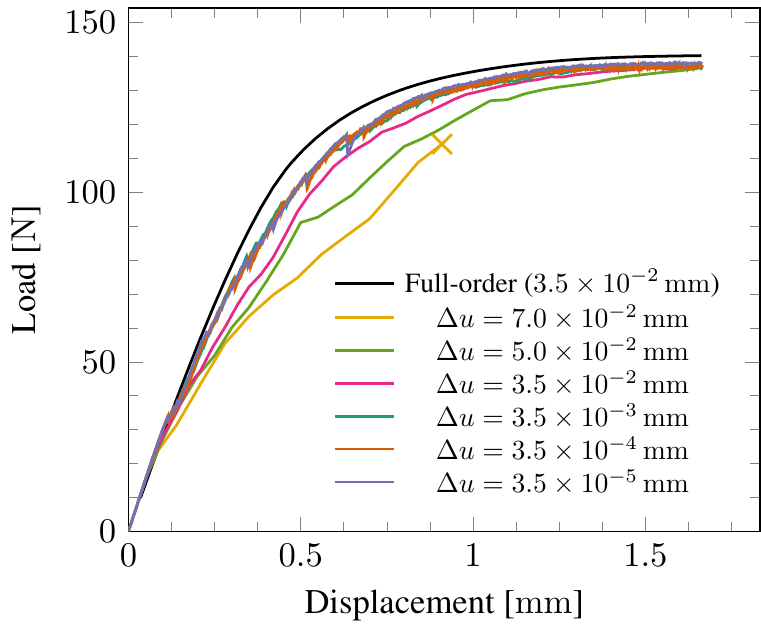}
  \caption{\fetwo\ example: Effect of time step size on surrogate predictions.}
  \label{FIholedogtimesteps}
\end{figure}

As one final result, we run the monotonic simulation with the hybrid surrogate for different time step sizes. As previously mentioned, the hybrid approach allows for explicit update of \props\ within an implicit simulation by obtaining the tangent stiffness matrix directly from the decoder. This however introduces a time step size dependency whose impact merits investigation. We plot in \cref{FIholedogtimesteps} predictions with step sizes spanning four orders of magnitude, including the same one used to obtain the full-order response. The combination of the explicit property update with the online stabilization procedure indeed introduces an upper bound for time step size for this specific problem. It stands to reason that the sensitivity to time step size also depends on the choice of decoder and on which material properties are included in \props. Further investigation into the matter in future works is therefore warranted.

\section{Conclusions}
\label{SEconclusions}

In this paper, we propose a hybrid surrogate modeling architecture for multiscale modeling of heterogeneous materials. The model is composed of a data-driven encoder for material properties and a physics-based decoder that computes stresses. In the resulting architecture, the encoder increases the flexibility of existing material models by letting their properties evolve in time, while the decoder provides beneficial bias and interpretability to the model. The model is conceived with flexibility in mind, allowing existing implementations of physics-based material models to be used with no extra modifications. Furthermore, by letting the decoder directly receive strain inputs, the encoder architecture is highly flexible and allows for preservation of frame independence. A semi-explicit online prediction algorithm is also proposed that allows for imposing extra constraints to model behavior in a semi-supervised way.

We demonstrate the architecture by reproducing pressure-dependent elastoplastic behavior coming from homogenized fiber-reinforced composite micromodels. The simple model with a linear-elastic decoder learned faster than conventional data-driven surrogates, allowed for lossless feature space dimensionality reduction through the use of strain invariants, and was able to approximate path-dependent behavior through a simple history-aware feature extractor. Models with perfectly-plastic $J_2$ decoders were shown to successfully learn nonlinear hardening and pressure dependency and predict unloading-reloading while being trained exclusively on monotonic data, outperforming a state-of-the-art mesomodel for composites in accuracy for arbitrary loading directions. Employing as decoder the same plasticity model used at the microscale led to highly-accurate monotonic response and fairly accurate extrapolation to unloading/reloading behavior. Finally, the model was used to solve a complex \fetwo\ model and the benefit of the online stabilization procedure was demonstrated.

We find the approach to be a promising new way to build hybrid surrogates which therefore merits further research on a number of fronts.
The current architecture is not by construction concerned with enforcing unconditional thermodynamic consistency or other physical constraints of interest. Although we do find empirically that well-trained surrogates with thermodynamically consistent decoders tend to perform well, some constitutive models might not be suitable for having their properties evolve in time. Fortunately, the framework can cope with extra constraints without necessarily giving up on its flexibility, by enforcing them locally through online retraining. Although training exclusively on monotonic paths already allows for path dependency to be fairly well captured, some decoders might perform better in extrapolation if trained with a (small) number of extra non-monotonic and non-proportional strain paths --- for instance when encoder and decoder can each explain the same phenomenon on their own (\eg pressure dependency in the model by Melro \textit{et al.}). We also foresee combining the present approach with the one in \cite{maiaPhysicallyRecurrentNeural2022} into a unified family of flexible hybrid surrogates with a range of possible combinations of feature extractors for physics-rich time convolution, fixed-property models with learned strain distributions and evolving material models. 



\section*{Acknowledgements}

The authors gratefully acknowledge the TU Delft AI Initiative for their support through the SLIMM AI Lab. FM also acknowledges financial support from the Netherlands Organization for Scientific Research (NWO) under Vidi grant nr. 16464.


\bibliographystyle{unsrt}
\bibliography{paper}

\end{document}